\documentclass[11pt]{amsart}

\usepackage{amssymb,amsmath,a4wide,graphicx,stmaryrd}
\usepackage[pagebackref=true]{hyperref}
\hypersetup{backref}
\usepackage[all]{xy}

\title{Directed algebraic topology and higher dimensional transition
  systems}

\author[P. Gaucher]{Philippe Gaucher}

\address{Laboratoire PPS  (CNRS UMR 7126)\\
  Universit{\'e} Paris 7--Denis Diderot\\
  Site Chevaleret\\ Case 7014\\ 75205 PARIS Cedex 13 \\ France}

\urladdr{http://www.pps.jussieu.fr/{\~{}}gaucher/} 



\swapnumbers

\newcommand{\C}{\mathcal{C}}
\newcommand{\D}{\mathcal{D}}

\newcommand{\N}{\mathbb{N}}

\newcommand{\T}{\mathbb{T}}
\newcommand{\de}{\partial}
\newcommand{\p}\times
\renewcommand{\vec}{\overrightarrow}
\renewcommand{\P}{\mathbb{P}}

\newtheorem*{thmN}{Theorem}

\newtheorem{thm}{Theorem}[section]
\newtheorem{prop}[thm]{Proposition}
\newtheorem{lem}[thm]{Lemma}
\newtheorem{cor}[thm]{Corollary}

\newtheorem{defn}[thm]{Definition}
\newtheorem{nota}[thm]{Notation}

\newcommand{\bd}{\begin{defn}}
\newcommand{\ed}{\end{defn}}
\newcommand{\bp}{\begin{prop}}
\newcommand{\ep}{\end{prop}}
\newcommand{\bth}{\begin{thm}}
\renewcommand{\eth}{\end{thm}}
\newcommand{\bpf}{\begin{proof}}
\newcommand{\epf}{\end{proof}}

\newcommand{\fl}[1]{\ar@{->}[ll]_-{#1}}
\newcommand{\fr}[1]{\ar@{->}[rr]^-{#1}}
\newcommand{\fd}[1]{\ar@{->}[dd]_-{#1}}
\newcommand{\fu}[1]{\ar@{->}[uu]^-{#1}}
\newcommand{\f}[2]{\ar@{->}[#1]|{#2}}
\newcommand{\ff}[2]{\ar@2{->}[#1]|{#2}}
\newcommand{\frr}[1]{\ar@{->}[rrrr]^-{#1}}

\renewcommand{\top}{{\mathbf{Top}}}

\newcommand{\iso}{\cong}
\newcommand{\lp}{\left(}
\newcommand{\rp}{\right)}
\newcommand{\ot}{\otimes}

\renewcommand{\leq}{\leqslant}
\renewcommand{\geq}{\geqslant}

\def\cartesien{%
  \ar@{-}[]+R+<6pt,-2pt>;[]+RD+<6pt,-6pt>%
  \ar@{-}[]+D+<2pt,-6pt>;[]+RD+<6pt,-6pt>%
}
\def\cocartesien{%
  \ar@{-}[]+L+<-6pt,+2pt>;[]+LU+<-6pt,+6pt>%
  \ar@{-}[]+U+<-2pt,+6pt>;[]+LU+<-6pt,+6pt>%
}
\def\hocartesien{%
  \ar@{-}[]+R+<6pt,-2pt>;[]+RD+<6pt,-6pt>_{h}%
  \ar@{-}[]+D+<2pt,-6pt>;[]+RD+<6pt,-6pt>%
}
\def\hococartesien{%
  \ar@{-}[]+L+<-6pt,+2pt>;[]+LU+<-6pt,+6pt>_{h}%
  \ar@{-}[]+U+<-2pt,+6pt>;[]+LU+<-6pt,+6pt>%
}

\newcommand{\brm}[1]{\rm{\mathbf{#1}}}

\newcommand{\dtop}{{\brm{Flow}}}
\newcommand{\set}{{\brm{Set}}}
\newcommand{\poset}{{\brm{PoSet}}}

\newcommand{\proc}{{\brm{Proc}}}

\newcommand{\hda}{{\brm{HDA}}}

\DeclareMathOperator{\rec}{rec}
\DeclareMathOperator{\id}{Id}
\DeclareMathOperator{\sh}{Sh}

\newcommand{\liminj}{\varinjlim}
\newcommand{\limproj}{\varprojlim}
\newcommand{\cat}{{\mathbf{Cat}}}

\DeclareMathOperator{\hdts}{\brm{HDTS}}
\DeclareMathOperator{\whdts}{\brm{WHDTS}}
\DeclareMathOperator{\Mod}{Mod}
\DeclareMathOperator{\cub}{\underline{Cub}}
\DeclareMathOperator{\cube}{CUBE}
\DeclareMathOperator{\opt}{Opt}

\makeatletter
\def\varholim@#1#2{%
  \vtop{\m@th\ialign{##\cr
    \hfil$#1\operator@font holim$\hfil\cr
    \noalign{\nointerlineskip\kern1.5\ex@}#2\cr
    \noalign{\nointerlineskip\kern-\ex@}\cr}}%
}
\def\holimproj{%
  \mathop{\mathpalette\varholim@{\leftarrowfill@\textstyle}}\nmlimits@
}
\def\holiminj{%
  \mathop{\mathpalette\varholim@{\rightarrowfill@\textstyle}}\nmlimits@
}
\makeatother

\DeclareMathOperator{\cosk}{cosk}
\DeclareMathOperator{\COSK}{\vec{\cosk}}

\DeclareMathOperator{\cell}{{\brm{cell}}}
\DeclareMathOperator{\cof}{{\brm{cof}}}
\DeclareMathOperator{\inj}{{\brm{inj}}}
\newcommand{\ddownarrow}{{\downarrow}}

\begin{document}

\begin{abstract}
  Cattani-Sassone's notion of higher dimensional transition system is
  interpreted as a small-orthogonality class of a locally finitely
  presentable topological category of weak higher dimensional
  transition systems. In particular, the higher dimensional transition
  system associated with the labelled $n$-cube turns out to be the
  free higher dimensional transition system generated by one
  $n$-dimensional transition. As a first application of this
  construction, it is proved that a localization of the category of
  higher dimensional transition systems is equivalent to a locally
  finitely presentable reflective full subcategory of the category of
  labelled symmetric precubical sets. A second application is to
  Milner's calculus of communicating systems (CCS): the mapping taking
  process names in CCS to flows is factorized through the category of
  higher dimensional transition systems. The method also applies to
  other process algebras and to topological models of concurrency
  other than flows.
\end{abstract}

\maketitle

\tableofcontents

\section{Introduction}

\subsection*{Presentation of the results}

In directed algebraic topology, the concurrent execution of $n$
actions is modelled by a full $n$-cube, each coordinate corresponding
to one of the $n$ actions. In this setting, a general concurrent
process is modelled by a gluing of $n$-cubes modelling the execution
paths and the higher dimensional homotopies between them. Various
topological models are being studied: in alphabetic and
non-chronological order, $d$-space \cite{mg}, $d$-space generated by
cubes \cite{FR}, flow \cite{model3}, globular complex \cite{diCW},
local po-space \cite{MR1683333}, locally preordered space \cite{SK},
multipointed $d$-space \cite{interpretation-glob}, and more
\cite{survol} (this list is probably not complete, indeed).  The
combinatorial model of labelled (symmetric) precubical set is also of
interest because, with such a model, it is exactly known where the
cubes are located in the geometry of the object. It was introduced for
the first time in \cite{labelled} \cite{exHDA}, following ideas from
\cite{EWDCooperating} \cite{Pratt} \cite{Gunawardena1} \cite{rvg} (the
last paper is a recent survey containing references to older papers),
and improved in \cite{ccsprecub} \cite{symcub} in relation with the
study of process algebras. The paper \cite{ccsprecub} treated the case
of labelled precubical sets, and the paper \cite{symcub} the more
general cases of labelled symmetric precubical sets and labelled
symmetric transverse precubical sets.

An apparently different philosophy is the one of higher dimensional
transition system. This notion, introduced in \cite{MR1461821}, models
the concurrent execution of $n$ actions by a transition between two
states labelled by a multiset of $n$ actions. A multiset is a set with
possible repetition of some elements (e.g. $\{0,0,2,3,3,3\}$). It is
usually modelled by an object of $\set\ddownarrow \N^*$, i.e. by a set
map $N : X\rightarrow \N^*$ where $X$ is the underlying set of the
multiset $N$ in which $x\in X$ appears $N(x)>0$ times. A higher
dimensional transition system must satisfy several natural axioms
CSA1, CSA2 and CSA3 (cf. Definition~\ref{def_hdts}). This notion is a
generalization of the $1$-dimensional notion of transition system in
which transitions between states are labelled by one action (e.g.,
\cite[Section~2.1]{MR1365754}). The latter $1$-dimensional notion
cannot of course model concurrency.

One of the purposes of this paper is to make precise the link between
process algebras modelled as labelled symmetric precubical sets, as
higher dimensional transition systems, and as flows, by introducing
the notion of weak higher dimensional transition system. The only
process algebras treated are the ones in Milner's calculus of
communicating systems (CCS) \cite{MR1365754} \cite{0683.68008}. And
only the topological model of flows introduced in \cite{model3} is
used. Similar results can easily be obtained for other process
algebras and for topological models of concurrency other than
flows. For other synchronization algebras, one needs only to change
the set of synchronizations in Definition~\ref{fibered}. For other
topological models of concurrency one needs only to change the
realization of the full $n$-cube $[n]^{cof}$ in
Definition~\ref{geom_rea}.  These modifications do not affect the
mathematical results of the paper. The first main result can then be
stated as follows:

\begin{thmN} (Theorem~\ref{colim1}, Theorem~\ref{colim2},
  Theorem~\ref{factor} and Corollary~\ref{final}) The mapping defined
  in \cite{ccsprecub} and \cite{symcub} taking each CCS process name
  to the geometric realization as flow $|\square_S\llbracket
  P\rrbracket|_{flow}$ of the labelled symmetric precubical set
  $\square_S\llbracket P\rrbracket$ factors through Cattani-Sassone's
  category of higher dimensional transition systems.
\end{thmN}

In fact, the functorial factorization $|\T(K)| \iso |K|_{flow}$ exists
as soon as $K$ satisfies the HDA paradigm and $\T(K)$ the Unique
intermediate state axiom (Every $K$ satisfying the latter condition is
called a strong labelled symmetric precubical set).

Let us recall for the reader that the semantics of process algebras
used in this paper in Section~\ref{proc} is the one of
\cite{symcub}. This semantics is nothing else but the labelled free
symmetric precubical set generated by the labelled precubical set
given in \cite{ccsprecub}. The reason for working with labelled
symmetric precubical sets in this paper is that this category is
closely related to the category of (weak) higher dimensional
transition systems by Theorem~\ref{ff}: the full subcategories in the
two categories generated by the labelled $n$-cubes for all $n\geq 0$
are isomorphic.

The interest of the combinatorial model of (weak) higher dimensional
transition systems is that the HDA paradigm
(cf. Section~\ref{paradigm}) is automatically satisfied.  That is to
say, the concurrent execution of $n$ actions (with $n\geq 2$) always
assembles to exactly one $n$-cube in a (weak) higher dimensional
transition system. Indeed, the realization functor $\T$ from labelled
symmetric precubical sets to weak higher dimensional transition
systems factors through the category of labelled symmetric precubical
sets satisfying the HDA paradigm by Theorem~\ref{colim2}.  On the
contrary, as already explained in \cite{ccsprecub} and in
\cite{symcub}, there exist labelled (symmetric) precubical sets
containing $n$-tuples of actions running concurrently which assemble
to several different $n$-cubes. Let us explain this phenomenon for the
case of the square. Consider the concurrent execution of two actions
$a$ and $b$ as depicted in Figure~\ref{concab}. Let $S = \{0,1\}\p
\{0,1\}$ be the set of states. Let $L = \{a,b\}$ be the set of actions
with $a \neq b$. The boundary of the square is modelled by adding to
the set of states $S$ the four $1$-transitions $((0,0),a,(1,0))$,
$((0,1),a,(1,1))$, $((0,0),b,(0,1))$ and $((1,0),b,(1,1))$. The
concurrent execution of $a$ and $b$ is modelled by adding the
$2$-transitions $((0,0),a,b,(1,1))$ and $((0,0),b,a,(1,1))$. Adding
one more time the two $2$-transitions $((0,0),a,b,(1,1))$ and
$((0,0),b,a,(1,1))$ does not change anything to the object since the
set of transitions remains equal to
\begin{multline*}\{((0,0),a,(1,0)),((0,1),a,(1,1)),((0,0),b,(0,1)),
  ((1,0),b,(1,1)),\\((0,0),a,b,(1,1)),((0,0),b,a,(1,1))\}.\end{multline*}
On the contrary, the labelled symmetric precubical set $\square_S[a,b]
\sqcup_{\de{\square_S}[a,b]} \square_S[a,b]$ contains two different
labelled squares $\square_S[a,b]$ modelling the concurrent execution
of $a$ and $b$, obtaining this way a geometric object homotopy
equivalent to a $2$-dimensional sphere (see
Proposition~\ref{ecrasement_invisible}). This is meaningless from a
computer scientific point of view. Indeed, either the two actions $a$
and $b$ run sequentially, and the square must remain empty, or the two
actions $a$ and $b$ run concurrently and the square must be filled by
exactly one square modelling concurrency. The topological hole created
by the presence of two squares as in $\square_S[a,b]
\sqcup_{\de{\square_S}[a,b]} \square_S[a,b]$ does not have any
computer scientific interpretation. The concurrent execution of two
actions (and more generally of $n$ actions) must be modelled by a
contractible object.

The factorization of $\T$ even yields a faithful functor
$\overline{\T}$ from labelled symmetric precubical sets satisfying the
HDA paradigm to weak higher dimensional transition systems by
Corollary~\ref{injection0}. However, the functor $\overline{\T}$ is
not full by Proposition~\ref{eq}. It only induces an equivalence of
categories by restricting to a full subcategory:

\begin{thmN} (Theorem~\ref{application1}) The localization of the
  category of higher dimensional transition systems by the
  cubification functor is equivalent to a locally finitely presentable
  reflective full subcategory of the category of labelled symmetric
  precubical sets. In this localization, two higher dimensional
  transition systems are isomorphic if they have the same cubes and
  they only differ by their set of actions.
\end{thmN}

We must introduce the technical notion of weak higher dimensional
transition system since there exist labelled symmetric precubical sets
$K$ such that $\T(K)$ is not a higher dimensional transition system by
Proposition~\ref{not-strong}. It is of course not difficult to find a
labelled symmetric precubical set contradicting CSA1 of
Definition~\ref{defCS} (e.g., Figure~\ref{Da}). It is also possible to
find counterexamples for the other axioms CSA2 and CSA3 of higher
dimensional transition system. This matters: if a labelled symmetric
precubical set $K$ is such that $\T(K)$ is not a higher dimensional
transition system, then it cannot be constructed from a process
algebra.

\subsection*{Organization of the paper}

Section~\ref{weakhdts} expounds the notion of weak higher dimensional
transition system. The notion of multiset recalled in the introduction
is replaced by the Multiset axiom on tuples to make the categorical
treatment easier. Logical tools are used to prove that the category of
weak higher dimensional transition systems is locally finitely
presentable and topological. Section~\ref{def_hdts} recalls
Cattani-Sassone's notion of higher dimensional transition system. It
is proved that every higher dimensional transition system is a weak
one. The notion of higher dimensional transition system is also
reformulated to make it easier to use. The Unique intermediate state
axiom is introduced for that purpose. It is also proved in the same
section that the set of transitions of any reasonable colimit is the
union of the transitions of the components
(Theorem~\ref{cube_final}). It is proved in Section~\ref{small_class}
that higher dimensional transition systems assemble to a
small-orthogonality class of the category of weak higher dimensional
transition systems (Corollary~\ref{ortho_hdts_class}). This implies
that the category of higher dimensional transition systems is a full
reflective locally finitely presentable category of the category of
weak higher dimensional transition
systems. Section~\ref{rappel_sym_pre_set} recalls the notion of
labelled symmetric precubical set. This section collects information
scattered between \cite{ccsprecub} and
\cite{symcub}. Section~\ref{paradigm} defines the paradigm of higher
dimensional automata (HDA paradigm). It is the adaptation to the
setting of labelled symmetric precubical sets of the analogous
definition presented in \cite{ccsprecub} for labelled precubical
sets. A labelled symmetric precubical set satisfies the HDA paradigm
if every labelled $p$-shell with $p\geq 1$ can be filled by at most
one labelled $(p+1)$-cube. This notion is a technical tool for various
proofs of this paper. It is proved in the same section that the full
subcategory of labelled symmetric precubical sets satisfying the HDA
paradigm is a full reflective subcategory of the category of labelled
symmetric precubical sets by proving that it is a small-orthogonality
class as well. It is also checked in the same section that the full
labelled $n$-cube satisfies the HDA paradigm (this trivial point is
fundamental!). Section~\ref{iso_cube} establishes that the full
subcategory of labelled $n$-cubes of the category of labelled
symmetric precubical sets is isomorphic to the full subcategory of
labelled $n$-cubes of the category of (weak) higher dimensional
transition systems (Theorem~\ref{ff}). The proof is of combinatorial
nature. Section~\ref{realization} constructs the realization functor
from labelled symmetric precubical sets to weak higher dimensional
transition systems. And it is proved that this functor factors through
the full subcategory of labelled symmetric precubical sets satisfying
the HDA paradigm. The two functors, the realization functor and its
factorization are left adjoints (Theorem~\ref{colim1} and
Theorem~\ref{colim2}). Section~\ref{property_func} studies when these
latter functors are faithful and full. It is proved that the HDA
paradigm is related to faithfulness and that the combination of the
HDA paradigm together and the Unique intermediate state axiom is
related to fullness.  Section~\ref{coreflec} uses all previous results
to compare the two settings of higher dimensional transition systems
and labelled symmetric precubical sets.  Section~\ref{geo} is a
straightforward but crucial application of the previous results. It is
proved in Theorem~\ref{factor} that the geometric realization as flow
of a labelled symmetric precubical set $K$ is the geometric
realization as flow of its realization as weak higher dimensional
transition system provided that $K$ is strong and satisfies the HDA
paradigm. The purpose of Section~\ref{proc} is to prove that these
conditions are satisfied by the labelled symmetric precubical sets
coming from process algebras. Hence we obtain the second application
stated in Corollary~\ref{final}.

\section{Prerequisites} 

The notations used in this paper are standard. A small class is called
a set.  All categories are locally small. The set of morphisms from
$X$ to $Y$ in a category $\C$ is denoted by $\C(X,Y)$. The identity of
$X$ is denoted by $\id_X$. Colimits are denoted by $\liminj$ and
limits by $\limproj$.

The reading of this paper requires general knowledge in category
theory \cite{MR1712872}, in particular in presheaf theory
\cite{MR1300636}, but also a good understanding of the theory of
locally presentable categories \cite{MR95j:18001} and of the theory of
topological categories \cite{topologicalcat}. A few model category
techniques are also used \cite{MR1361887} \cite{MR99h:55031}
\cite{ref_model2} in the proof of Theorem~\ref{quotient} and in
Section~\ref{geo}.

A short introduction to process algebra can be found in
\cite{MR1365754}. An introduction to CCS (Milner's calculus of
communicating systems \cite{0683.68008}) for mathematician is
available in \cite{ccsprecub} and in \cite{symcub}. Hardly any
knowledge of process algebra is needed to read Section~\ref{proc} of
the paper. In fact, the paper \cite{ccsprecub} can be taken as a
starting point.

Some salient mathematical facts are collected in this section. Of
course, this section does not intend to be an introduction to these
notions. It will only help the reader to understand what kinds of
mathematical tools are used in this work.

Let $\lambda$ be a regular cardinal (see for example \cite[p
160]{MR1697766}). When $\lambda=\aleph_0$, the word ``$\lambda-$'' is
replaced by the word ``finitely''. An object $C$ of a category $\C$ is
$\lambda$-presentable when the functor $\C(C,-)$ preserves
$\lambda$-directed colimits. Practically, that means that every map $C
\rightarrow \liminj C_i$ factors as a composite $C \rightarrow C_i
\rightarrow \liminj C_i$ when the colimit is $\lambda$-directed.  A
$\lambda$-accessible category is a category having $\lambda$-directed
colimits such that each object is generated (in some strong sense) by
a set of $\lambda$-presentable objects. For example, each object is a
$\lambda$-directed colimit of a subset of a given set of
$\lambda$-presentable objects. If moreover the category is cocomplete,
it is called a locally $\lambda$-presentable category. We use at
several places of the paper a logical characterization of accessible
and locally presentable categories which are axiomatized by theories
with set of sorts $\{s\} \cup \Sigma$, $s$ being the sort of states
and $\Sigma$ a non-empty fixed set of labels. Another kind of locally
presentable category is a category of presheaves, and any comma
category constructed from it.  Every locally presentable category has
a set of generators, is complete, cocomplete, wellpowered and
co-wellpowered. The Special Adjoint Functor Theorem SAFT is then
usable to prove the existence of right adjoints. A functor between
locally $\lambda$-presentable category is $\lambda$-accessible if it
preserves $\lambda$-directed colimits (or equivalently
$\lambda$-filtered colimits). Another important fact is that a functor
between locally presentable categories is a right adjoint if and only
if it is accessible and limit-preserving.

An object $C$ is orthogonal to a map $X\rightarrow Y$ if every map
$X\rightarrow C$ factors uniquely as a composite $X \rightarrow Y
\rightarrow C$. A full subcategory of a given category is reflective
if the inclusion functor is a right adjoint. The left adjoint to the
inclusion is called the reflection. In a locally presentable category,
the full subcategory of objects orthogonal to a given set of morphisms
is an example of a reflective subcategory. Such a category, called a
small-orthogonality class, is locally presentable. And the inclusion
functor is of course accessible and limit-preserving.

The paradigm of topological category over the category of $\set$ is
the one of general topological spaces with the notions of initial
topology and final topology. More precisely, a functor $\omega: \C
\rightarrow \D$ is topological if each cone $(f_i:X \rightarrow \omega
A_i)_{i\in I}$ where $I$ is a class has a unique $\omega$-initial lift
(the initial structure) $(\overline{f}_i:A \rightarrow A_i)_{i\in I}$,
i.e.: 1) $\omega A= X$ and $\omega \overline{f}_i = f_i$ for each
$i\in I$; 2) given $h: \omega B \rightarrow X$ with $f_i h= \omega
\overline{h}_i$, $\overline{h}_i:B\rightarrow A_i$ for each $i\in I$,
then $h=\omega \overline{h}$ for a unique $\overline{h} : B
\rightarrow A$.  Topological functors can be characterized as functors
such that each cocone $(f_i:\omega A_i \rightarrow X)_{i\in I}$ where
$I$ is a class has a unique $\omega$-final lift (the final structure)
$\overline{f}_i:A_i \rightarrow A$, i.e.: 1) $\omega A= X$ and $\omega
\overline{f}_i = f_i$ for each $i\in I$; 2) given $h: X \rightarrow
\omega B$ with $hf_i = \omega \overline{h}_i$, $\overline{h}_i:A_i
\rightarrow B$ for each $i\in I$, then $h=\omega \overline{h}$ for a
unique $\overline{h} : A \rightarrow B$. Let us suppose $\D$ complete
and cocomplete. A limit (resp. colimit) in $\C$ is calculated by
taking the limit (resp. colimit) in $\D$, and by endowing it with the
initial (resp. final) structure. In this work, a topological category
is a topological category over the category $\set^{\{s\} \cup \Sigma}$
where $\{s\} \cup \Sigma$ is as above the set of sorts.

Let $i:A\longrightarrow B$ and $p:X\longrightarrow Y$ be maps in a
category $\C$. Then $i$ has the left lifting property (LLP)
with respect to $p$ (or $p$ has the right lifting property
(RLP) with respect to $i$) if for every commutative square
\[
\xymatrix{
A\fd{i} \fr{\alpha} && X \fd{p} \\
&&\\
B \ar@{-->}[rruu]^{g}\fr{\beta} && Y}
\]
there exists a lift $g$ making both triangles commutative.

Let $\C$ be a cocomplete category.  If $K$ is a set of morphisms of
$\C$, then the class of morphisms of $\C$ that satisfy the RLP with
respect to every morphism of $K$ is denoted by $\inj(K)$ and the class
of morphisms of $\C$ that are transfinite compositions of pushouts of
elements of $K$ is denoted by $\cell(K)$. Denote by $\cof(K)$ the
class of morphisms of $\C$ that satisfy the LLP with respect to the
morphisms of $\inj(K)$.  It is a purely categorical fact that
$\cell(K)\subset \cof(K)$. Moreover, every morphism of $\cof(K)$ is a
retract of a morphism of $\cell(K)$ as soon as the domains of $K$ are
small relative to $\cell(K)$ \cite[Corollary~2.1.15]{MR99h:55031}. An
element of $\cell(K)$ is called a relative $K$-cell complex.  If $X$
is an object of $\C$, and if the canonical morphism
$\varnothing\longrightarrow X$ is a relative $K$-cell complex, then
the object $X$ is called a $K$-cell complex.

Let $\C$ be a category. A weak factorization system is a pair
$(\mathcal{L},\mathcal{R})$ of classes of morphisms of $\C$ such that
the class $\mathcal{L}$ is the class of morphisms having the LLP with
respect to $\mathcal{R}$, such that the class $\mathcal{R}$ is the
class of morphisms having the RLP with respect to $\mathcal{L}$ and
such that every morphism of $\C$ factors as a composite $r\circ \ell$
with $\ell\in \mathcal{L}$ and $r\in \mathcal{R}$. The weak
factorization system is functorial if the factorization $r\circ \ell$
is a functorial factorization. It is cofibrantly generated if it is of
the form $(\cof(K),\inj(K))$ for some set of maps $K$.

A model category is a complete cocomplete category equipped with a
model structure consisting of three classes of morphisms ${\rm Cof}$,
${\rm Fib}$, $\mathcal{W}$ respectively called cofibration, fibration
and weak equivalence such that the pairs of classes of morphisms
$({\rm Cof},{\rm Fib }\cap \mathcal{W})$ and $({\rm Cof}\cap
\mathcal{W},{\rm Fib })$ are weak factorization systems and such that
if two of the three morphisms $f,g,g\circ f$ are weak equivalences,
then so is the third one. This model structure is cofibrantly
generated provided that the two weak factorization systems $({\rm
  Cof},{\rm Fib }\cap \mathcal{W})$ and $({\rm Cof}\cap
\mathcal{W},{\rm Fib })$ are cofibrantly generated. The only model
category used in this paper is the one of flows. We need only in fact
the notion of cofibrant replacement. For an object $X$ of a model
category, the canonical map $\varnothing \rightarrow X$ factors as a
composite $\varnothing \rightarrow X^{cof} \rightarrow X$ where the
left-hand map is a cofibration and the right-hand map a trivial
fibration, i.e. an element of ${\rm Fib }\cap \mathcal{W}$. The object
$X^{cof}$ is called a cofibrant replacement of $X$.

The proof of Theorem~\ref{quotient} uses the fact that for every set of
morphisms $K$ in a locally presentable category, a map $X \rightarrow 
Y$ always factors as a composite  $X \rightarrow Z \rightarrow
Y$ where the left-hand map is an object of $\cell(K)$ and the
right-hand map an object of $\inj(K)$. 

Beware of the fact that the word ``model'' has three different
meanings in this paper, a logical one, a homotopical one, and also a
non-mathematical one like in the sentence ``the $n$-cube models the
concurrent execution of $n$ actions''.

\section{Weak higher dimensional transition systems}
\label{weakhdts}

The formalism of multiset as used in \cite{MR1461821} is not easy to
handle. In this paper, an $n$-transition between two states $\alpha$
and $\beta$ (or from $\alpha$ to $\beta$) modelling the concurrent
execution of $n$ actions $u_1, \dots, u_n$ with $n\geq 1$ is modelled
by an $(n+2)$-tuple $(\alpha,u_1,\dots,u_n,\beta)$ satisfying the new
Multiset axiom: for every permutation $\sigma$ of $\{1,\dots,n\}$,
$(\alpha, u_{\sigma(1)}, \dots, u_{\sigma(n)}, \beta)$ is an
$n$-transition.

\begin{nota} We fix a nonempty set of {\rm labels} $\Sigma$. We
  suppose that $\Sigma$ always contains a distinguished element
  denoted by $\tau$. \end{nota}

\bd A {\rm weak higher dimensional transition system} consists of a
triple \[(S,\mu:L\rightarrow \Sigma,T=\bigcup_{n\geq 1}T_n)\] where
$S$ is a set of {\rm states}, where $L$ is a set of {\rm actions},
where $\mu:L\rightarrow \Sigma$ is a set map called the {\rm labelling
  map}, and finally where $T_n\subset S\p L^n\p S$ for $n \geq 1$ is a
set of {\rm $n$-transitions} or {\rm $n$-dimensional transitions} such
that one has:
\begin{itemize}
\item (Multiset axiom) For every permutation $\sigma$ of
  $\{1,\dots,n\}$ with $n\geq 2$, if $(\alpha,u_1,\dots,u_n,\beta)$ is
  a transition, then
  $(\alpha,u_{\sigma(1)},\dots,u_{\sigma(n)},\beta)$ is a transition
  as well.
\item (Coherence axiom) For every $(n+2)$-tuple
  $(\alpha,u_1,\dots,u_n,\beta)$ with $n\geq 3$, for every $p,q\geq 1$
  with $p+q<n$, if the five tuples $(\alpha,u_1,\dots,u_n,\beta)$,
  $(\alpha,u_1,\dots,u_p,\nu_1)$, $(\nu_1,u_{p+1},\dots,u_n,\beta)$,
  $(\alpha,u_1,\dots,u_{p+q},\nu_2)$ and
  $(\nu_2,u_{p+q+1},\dots,u_n,\beta)$ are transitions, then the
  $(q+2)$-tuple $(\nu_1,u_{p+1},\dots,u_{p+q},\nu_2)$ is a transition
  as well.
\end{itemize}
A map of weak higher dimensional transition systems
\[f:(S,\mu : L \rightarrow \Sigma,(T_n)_{n\geq 1}) \rightarrow
(S',\mu' : L' \rightarrow \Sigma ,(T'_n)_{n\geq 1})\] consists of a
set map $f_0: S \rightarrow S'$, a commutative square
\[
\xymatrix{
  L \ar@{->}[r]^-{\mu} \ar@{->}[d]_-{\widetilde{f}}& \Sigma \ar@{=}[d]\\
  L' \ar@{->}[r]_-{\mu'} & \Sigma}
\] 
such that if $(\alpha,u_1,\dots,u_n,\beta)$ is a transition, then
$(f_0(\alpha),\widetilde{f}(u_1),\dots,\widetilde{f}(u_n),f_0(\beta))$
is a transition. The corresponding category is denoted by $\whdts$.
The $n$-transition $(\alpha,u_1,\dots,u_n,\beta)$ is also called a
{\rm transition from $\alpha$ to $\beta$}.  \ed

\begin{nota} A transition $(\alpha,u_1,\dots,u_n,\beta)$ will be also
  denoted by $\alpha \stackrel{u_1,\dots,u_n}\longrightarrow
  \beta$. \end{nota}

\bth \label{loctop} The category $\whdts$ is locally finitely
presentable. The functor \[\omega : \whdts \longrightarrow
\set^{\{s\}\cup \Sigma}\] taking the weak higher dimensional
transition system $(S,\mu : L \rightarrow \Sigma,(T_n)_{n\geq 1})$ to
the $(\{s\}\cup \Sigma)$-tuple of sets $(S,(\mu^{-1}(x))_{x\in
  \Sigma}) \in \set^{\{s\}\cup \Sigma}$ is topological. \eth

\bpf Let $(f_i:\omega X_i \rightarrow (S,(L_x)_{x\in \Sigma}))_{i\in
  I}$ be a cocone where $I$ is a class with $X_i = (S_i,\mu_i:L_i
\rightarrow \Sigma, T^i = \bigcup_{n\geq 1} T^i_n)$. The closure by
the Multiset axiom and the Coherence axiom of the union of the images
of the $T^i$ in $\bigcup_{n\geq 1} (S \p L^n\p S)$ with $L =
\bigsqcup_{x\in \Sigma} L_x$ gives the final structure. Hence, the
functor $\omega$ is topological.

We use the terminology of \cite[Chapter~5]{MR95j:18001}.  Let us
consider the theory $\mathcal{T}$ in finitary first-order logic
defined by the set of sorts $\{s\}\cup \Sigma$, by a relational symbol
$T_{x_1,\dots,x_n}$ of arity $s\p x_1\p\dots\p x_n\p s$ for every
$n\geq 1$ and every $(x_1,\dots,x_n)\in \Sigma^n$, and by the axioms:
\begin{itemize}
\item for all $x_1,\dots,x_n\in \Sigma$, for all $n\geq 2$ and for all
  permutations $\sigma$ of $\{1,\dots,n\}$:
\[(\forall \alpha,u_1\dots,u_n,\beta),
T_{x_1,\dots,x_n}(\alpha,u_1,\dots,u_n,\beta) \Rightarrow
T_{x_{\sigma(1)},\dots,x_{\sigma(n)}}(\alpha,u_{\sigma(1)},\dots,u_{\sigma(n)},\beta).\]
\item for all $x_1,\dots,x_n\in \Sigma$, for all $n\geq 3$, for all
  $p,q\geq 1$ with $p+q<n$, 
\begin{multline*}(\forall \alpha,u_1\dots,u_n,\beta,\nu_1,\nu_2)
(T_{x_1,\dots,x_n}(\alpha,u_1,\dots,u_n,\beta) \wedge 
T_{x_1,\dots,x_p}(\alpha,u_1,\dots,u_p,\nu_1) \\\wedge 
T_{x_{p+1},\dots,x_n}(\nu_1,u_{p+1},\dots,u_n,\beta) \wedge 
T_{x_1,\dots,x_{p+q}}(\alpha,u_1,\dots,u_{p+q},\nu_2) \\\wedge
T_{x_{p+q+1},\dots,x_n}(\nu_2,u_{p+q+1},\dots,u_n,\beta)) \Rightarrow 
T_{x_{p+1},\dots,x_{p+q}}(\nu_1,u_{p+1},\dots,u_{p+q},\nu_2).\end{multline*} 
\end{itemize}
Since the axioms are of the form $(\forall {x}), \phi({x}) \Rightarrow
(\exists!  {y}\ \psi({x},{y}))$ (with no $y$) where $\phi$ and $\psi$
are conjunctions of atomic formulas with a finite number of arguments,
the category $\Mod(\mathcal{T})$ of models of $\mathcal{T}$ in
$\set^{\{s\}\cup \Sigma}$ is locally finitely presentable by
\cite[Theorem~5.30]{MR95j:18001}. It remains to observe that there is
an isomorphism of categories $\Mod(\mathcal{T}) \iso \whdts$ to
complete the proof.  \epf

Note that the category $\whdts$ is axiomatized by a universal strict
Horn theory without equality, i.e. by statements of the form $(\forall
{x}), \phi({x}) \Rightarrow \psi({x})$ where $\phi$ and $\psi$ are
conjunctions of atomic formulas without equalities. So
\cite[Theorem~5.3]{MR629337} provides another argument to prove that
the functor $\omega$ above is topological.

Let us conclude this section by some additional comments about
colimits in $\whdts$. We will come back to this question in
Theorem~\ref{cube_final}.

\bp \label{calcul_colim} Let $X=\liminj X_i$ be a colimit of weak
higher dimensional transition systems with $X_i = (S_i,\mu_i:L_i
\rightarrow \Sigma, T^i = \bigcup_{n\geq 1} T^i_n)$ and $X = (S,\mu:L
\rightarrow \Sigma, T = \bigcup_{n\geq 1} T_n)$. Then:
\begin{enumerate}
\item $S = \liminj S_i$, $L = \liminj L_i$, $\mu = \liminj \mu_i$
\item the union $\bigcup_{i} T^i$ of the image of the $T^i$ in
  $\bigcup_{n\geq 1} (S \p L^n\p S)$ satisfies the Multiset axiom.
\item $T$ is the closure of $\bigcup_{i} T^i$ under the Coherence
  axiom. 
\item when the union $\bigcup_{i} T^i$ is already closed under the
  Coherence axiom, this union is the final structure.
\end{enumerate}
\ep 

\bpf By \cite[Proposition~21.15]{topologicalcat}, (1) is a consequence
of the fact that $\whdts$ is topological over $\set^{\{s\}\cup
  \Sigma}$. (2) comes from the fact that each $T_i$ satisfies the
Multiset axiom. (4) is a consequence of (2). It remains to prove
(3). Let $G_0(\bigcup_{i} T^i) = \bigcup_{i} T^i$. Let us define
$G_\alpha(\bigcup_{i} T^i)$ by induction on the transfinite ordinal
$\alpha\geq 0$ by $G_\alpha(\bigcup_{i} T^i)= \bigcup_{\beta<\alpha}
G_\beta(\bigcup_{i} T^i)$ for every limit ordinal $\alpha$ and
$G_{\alpha+1}(\bigcup_{i} T^i)$ is obtained from $G_\alpha(\bigcup_{i}
T^i)$ by adding to $G_\alpha(\bigcup_{i} T^i)$ all $(q+2)$-tuples
$(\nu_1,u_{p+1},\dots,u_{p+q},\nu_2)$ such that there exist five
tuples $(\alpha,u_1,\dots,u_n,\beta)$, $(\alpha,u_1,\dots,u_p,\nu_1)$,
$(\nu_1,u_{p+1},\dots,u_n,\beta)$, $(\alpha,u_1,\dots,u_{p+q},\nu_2)$
and $(\nu_2,u_{p+q+1},\dots,u_n,\beta)$ of the set
$G_\alpha(\bigcup_{i} T^i)$. Hence we have the inclusions
$G_\alpha(\bigcup_{i} T^i) \subset G_{\alpha+1}(\bigcup_{i} T^i)
\subset \bigcup_{n\geq 1} (S \p L^n\p S)$ for all $\alpha\geq 0$. For
cardinality reason, there exists an ordinal $\alpha_0$ such that for
every $\alpha\geq \alpha_0$, one has $G_\alpha(\bigcup_{i} T^i) =
G_{\alpha_0}(\bigcup_{i} T^i)$.  By transfinite induction on
$\alpha\geq 0$, one sees that $G_\alpha(\bigcup_{i} T^i)$ satisfies
the Multiset axiom. So the closure $G_{\alpha_0}(\bigcup_{i} T^i)$ of
$\bigcup_{i} T^i$ under the Coherence axiom is the final structure and
$G_{\alpha_0}(\bigcup_{i} T^i) = T$.  \epf

\section{Higher dimensional transition systems}
\label{def_hdts}

Let us now propose our (slightly revised) version of higher
dimensional transition system.

\bd \label{defCS} A {\rm higher dimensional transition system} is a
triple \[(S,\mu:L\rightarrow \Sigma,T=\bigcup_{n\geq 1}T_n)\] where
$S$ is a set of {\rm states}, where $L$ is a set of {\rm actions},
where $\mu:L\rightarrow \Sigma$ is a set map called the {\rm labelling
  map}, and finally where $T_n\subset S\p L^n\p S$ is a set of {\rm
  $n$-transitions} or {\rm $n$-dimensional transitions} such that one
has:
\begin{enumerate}
\item (Multiset axiom) For every permutation $\sigma$ of
  $\{1,\dots,n\}$ with $n\geq 2$, if $(\alpha,u_1,\dots,u_n,\beta)$ is
  a transition, then $(\alpha, u_{\sigma(1)}, \dots, u_{\sigma(n)},
  \beta)$ is a transition as well.
\item (First Cattani-Sassone axiom CSA1) If $(\alpha,u,\beta)$ and
  $(\alpha,u',\beta)$ are two transitions such that $\mu(u) =
  \mu(u')$, then $u=u'$.
\item (Second Cattani-Sassone axiom CSA2) For every $n\geq 2$, every
  $p$ with $1\leq p<n$, and every transition $(\alpha, u_1, \dots,
  u_n, \beta)$, there exists a unique state $\nu_1$ and a unique state
  $\nu_2$ such that $(\alpha, u_{1}, \dots, u_{p}, \nu_1)$, $(\nu_1,
  u_{p+1}, \dots, u_{n}, \beta)$, $(\alpha, u_{p+1}, \dots, u_{n},
  \nu_2)$ and $(\nu_2, u_{1}, \dots, u_{p}, \beta)$ are transitions.
\item (Third Cattani-Sassone axiom CSA3) For every state $\alpha,
  \beta, \nu_1,\nu_2,\nu'_1,\nu'_2$ and every action $u_1,\dots,u_n$,
  with $p,q \geq 1$ and $p+q<n$, if the nine tuples
\begin{multline*}
  (\alpha,u_1,\dots,u_n,\beta), (\alpha,u_1,\dots,u_p,\nu_1),
  (\nu_1,u_{p+1},\dots,u_n,\beta),\\
  (\nu_1,u_{p+1},\dots,u_{p+q},\nu_2),
  (\nu_2,u_{p+q+1},\dots,u_n,\beta), (\alpha,u_{1},\dots,u_{p+q},\nu'_2),\\
  (\nu'_2,u_{p+q+1},\dots,u_n,\beta),
  (\alpha,u_{1},\dots,u_p,\nu'_1),
  (\nu'_1,u_{p+1},\dots,u_{p+q},\nu'_2)
\end{multline*}
are transitions, then $\nu_1 = \nu'_1$ and $\nu_2 = \nu'_2$.
\end{enumerate}
\ed 

Note that our notion of morphism of higher dimensional transition
systems differs from Cattani-Sassone's one: we take only the morphisms
between the underlying sets of states and actions preserving the
structure. This is necessary to develop the theory presented in this
paper. So it becomes false that two general higher dimensional
transition systems differing only by the set of actions are
isomorphic. However, this latter fact is true in some appropriate
categorical localization (see the very end of
Section~\ref{coreflec}). We also have something similar for (weak)
higher dimensional transition systems coming from strong labelled
symmetric precubical sets by Corollary~\ref{label_differ}, that is to
say from any labelled symmetric precubical set coming from process
algebras by Theorem~\ref{strong}.

Let us cite \cite{MR1461821}: ``CSA1 in the above definition simply
guarantees that there are no two transitions between the same states
carrying the same multiset of labels. CSA2 guarantees that all the
interleaving of a transition
$\alpha\stackrel{u_1,\dots,u_n}\longrightarrow \beta$ are present as
paths from $\alpha$ to $\beta$, whilst CSA3 ensures that such paths
glue together properly: it corresponds to the cubical laws of higher
dimensional automata''.

\bp \label{inclusion_weak} Every higher dimensional transition system
is a weak higher dimensional transition system. \ep

\bpf Let $X = (S,\mu:L\rightarrow \Sigma,T=\bigcup_{n \geq 1}T_n)$ be
a higher dimensional transition system. Let $(\alpha, u_1, \dots, u_n,
\beta)$ be a transition with $n \geq 3$. Let $p,q \geq 1$ with $p+q <
n$. Suppose that the five tuples $(\alpha,u_1,\dots,u_n,\beta)$,
$(\alpha,u_1,\dots,u_p,\nu_1)$, $(\nu_1,u_{p+1},\dots,u_n,\beta)$,
$(\alpha,u_1,\dots,u_{p+q},\nu_2)$ and
$(\nu_2,u_{p+q+1},\dots,u_n,\beta)$ are transitions. Let $\nu'_1$ be
the (unique) state of $X$ such that $(\alpha,u_1,\dots,u_{p},\nu'_1)$
and $(\nu'_1,u_{p+1},\dots,u_{p+q},\nu_2)$ are transitions of $X$.
Let $\nu'_2$ be the (unique) state of $X$ such that
$(\nu_1,u_{p+1},\dots,u_{p+q},\nu'_2)$ and
$(\nu'_2,u_{p+q+1},\dots,u_n,\beta)$ are transitions of $X$. Then
$\nu_1 = \nu'_1$ and $\nu_2 = \nu'_2$ by CSA3. Therefore the Coherence
axiom is satisfied.  \epf

\begin{nota} The full subcategory of higher dimensional transition
  systems is denoted by $\hdts$. So one has the inclusion $\hdts
  \subset \whdts$. \end{nota}

\bd A weak higher dimensional transition system satisfies the {\rm
  Unique Intermediate state axiom} if for every $n\geq 2$, every $p$
with $1\leq p<n$ and every transition $(\alpha,u_1,\dots,u_n,\beta)$,
there exists a unique state $\nu$ such that both the tuples
$(\alpha,u_1,\dots,u_p,\nu)$ and $(\nu,u_{p+1},\dots,u_n,\beta)$ are
transitions.  \ed

\bp A weak higher dimensional transition system satisfies the second
and third Catttani-Sassone axioms if and only if it satisfies the
Unique intermediate state axiom. \ep

\bpf A weak higher dimensional transition system satisfying CSA2 and
CSA3 clearly satisfies the Unique intermediate state axiom.
Conversely, if a weak higher dimensional transition system satisfies
the Unique intermediate state axiom, it clearly satisfies CSA2. Let
$\alpha, \beta, \nu_1,\nu_2,\nu'_1,\nu'_2$ be states and let
$u_1,\dots,u_n$ be actions with $n\geq 3$. Let $p,q\geq 1$ with
$p+q<n$. Suppose that
\begin{multline*}
  (\alpha,u_1,\dots,u_n,\beta), (\alpha,u_1,\dots,u_p,\nu_1),
  (\nu_1,u_{p+1},\dots,u_n,\beta),\\
  (\nu_1,u_{p+1},\dots,u_{p+q},\nu_2),
  (\nu_2,u_{p+q+1},\dots,u_n,\beta), (\alpha,u_{1},\dots,u_{p+q},\nu'_2),\\
  (\nu'_2,u_{p+q+1},\dots,u_n,\beta), (\alpha,u_{1},\dots,u_p,\nu'_1),
  (\nu'_1,u_{p+1},\dots,u_{p+q},\nu'_2)
\end{multline*}
are transitions. By the Coherence axiom, the tuple
$(\nu_1,u_{p+1},\dots,u_{p+q},\nu'_2)$ is a transition. By the Unique
intermediate state axiom, one obtains $\nu_1 = \nu'_1$ and $\nu_2 =
\nu'_2$. So CSA3 is satisfied too.  \epf

One obtains a new formulation of the notion of higher dimensional
transition system:

\bp \label{eq_hdts} A higher dimensional transition system is a weak
higher dimensional transition system satisfying CSA1 and the Unique
intermediate state axiom.  \ep

Let us conclude this section by an important remark about colimits of
weak higher dimensional transition systems satisfying the Unique
intermediate state axiom, so in particular about colimits of higher
dimensional transition systems.

\bth \label{cube_final} Let $X=\liminj X_i$ be a colimit of weak
higher dimensional transition systems such that every weak higher
dimensional transition system $X_i$ satisfies the Unique intermediate
state axiom. Let $X_i = (S_i,\mu_i:L_i \rightarrow \Sigma, T^i =
\bigcup_{n\geq 1} T^i_n)$ and $X = (S,\mu:L \rightarrow \Sigma, T =
\bigcup_{n\geq 1} T_n)$. Denote by $\bigcup_i T^i$ the union of the
images by the map $X_i\rightarrow X$ of the sets of transitions of the
$X_i$ for $i$ running over the set of objects of the base category of
the diagram $i\mapsto X_i$. Then the following conditions are
equivalent:
\begin{enumerate}
\item the weak higher dimensional transition system $X$ satisfies the
  Unique intermediate state axiom
\item the set of transitions $\bigcup_i T^i$ satisfies the Unique
  intermediate state axiom
\item the set of transitions $\bigcup_i T^i$ satisfies the Multiset
  axiom, the Coherence axiom and the Unique intermediate state axiom.
\end{enumerate}
Whenever one of the preceding three conditions is satisfied, the set
of transitions $\bigcup_i T^i$ is the final structure.  \eth

\bpf The set of transitions $\bigcup_i T^i$ always satisfies the
Multiset axiom by Proposition~\ref{calcul_colim}.

$(1)\Rightarrow (2)$. The set of transitions of $X$ is the closure
under the Coherence axiom of $\bigcup_i T^i$ by
Proposition~\ref{calcul_colim}. So $\bigcup_i T^i \subset T$.

$(2)\Rightarrow (3)$. Let $n\geq 3$. Let
$(\alpha,u_1,\dots,u_n,\beta)$ be a transition of $\bigcup_i T^i$. Let
$p,q\geq 1$ with $p+q<n$. Let $(\alpha,u_1,\dots,u_p,\nu_1)$,
$(\nu_1,u_{p+1},\dots,u_n,\beta)$, $(\alpha,u_1,\dots,u_{p+q},\nu_2)$
and $(\nu_2,u_{p+q+1},\dots,\linebreak[4]u_n,\beta)$ be four
transitions of $\bigcup_i T^i$. Let $i$ such that there exists a
transition $(\alpha^i,u^i_1,\dots,u^i_n,\beta^i)$ of $X_i$ taken by
the canonical map $X_i \rightarrow X$ to
$(\alpha,u_1,\dots,u_n,\beta)$. Since $X_i$ satisfies the Unique
intermediate state axiom, there exists a (unique) state $\nu_1^i$ and
a (unique) state $\nu_2^i$ of $X_i$ such that
$(\alpha^i,u^i_1,\dots,u^i_p,\nu^i_1)$,
$(\nu^i_1,u^i_{p+1},\dots,u^i_n,\beta^i)$,
$(\alpha^i,u^i_1,\dots,u^i_{p+q},\nu^i_2)$ and
$(\nu^i_2,u^i_{p+q+1},\dots,u^i_n,\beta)$ are four transitions of
$X_i$. Since $\bigcup_i T^i$ satisfies the Unique intermediate state
axiom as well, the map $X_i \rightarrow X$ takes $\nu^i_1$ to $\nu_1$
and $\nu^i_2$ to $\nu_2$. By the Coherence axiom applied to $X_i$, the
tuple $(\nu^i_1, u^i_{p+1},\dots,u^i_{p+q},\nu^i_2)$ is a transition
of $X_i$.  So the union $\bigcup_i T^i$ is closed under the Coherence
axiom. 

$(3)\Rightarrow (1)$. If $(3)$ holds, then the inclusion $\bigcup_i
T^i \subset T$ is an equality by
Proposition~\ref{calcul_colim}. Therefore the weak higher dimensional
transition system $X$ satisfies the Unique intermediate state axiom.

The last assertion is then clear. \epf

\section{Higher dimensional transition systems as a small-orthogonality
  class}
\label{small_class}

\begin{nota} Let $[0] = \{()\}$ and $[n] = \{0,1\}^n$ for $n \geq 1$.
  By convention, one has $\{0,1\}^0 = [0] = \{()\}$. The set $[n]$ is
  equipped with the product ordering $\{0<1\}^n$. \end{nota}

Let us now describe the higher dimensional transition system
associated with the $n$-cube for $n\geq 0$. 

\bp \label{cas_cube} Let $n\geq 0$ and $a_1,\dots,a_n\in \Sigma$. Let
$T_d\subset \{0,1\}^n \p \{(a_1,1),\dots,(a_n,n)\}^d \p \{0,1\}^n$
(with $d\geq 1$) be the subset of $(d+2)$-tuples
\[((\epsilon_1,\dots,\epsilon_n), (a_{i_1},i_1),\dots,(a_{i_d},i_d),
(\epsilon'_1,\dots,\epsilon'_n))\] such that
\begin{itemize}
\item $i_m = i_n$ implies $m = n$, i.e. there are no repetitions in the
  list $(a_{i_1},i_1),\dots,(a_{i_d},i_d)$
\item for all $i$, $\epsilon_i\leq \epsilon'_i$
\item $\epsilon_i\neq \epsilon'_i$ if and only if
  $i\in\{i_1,\dots,i_d\}$. 
\end{itemize}
Let $\mu : \{(a_1,1),\dots,(a_n,n)\} \rightarrow \Sigma$ be the set
map defined by $\mu(a_i,i) = a_i$. Then \[C_n[a_1,\dots,a_n] =
(\{0,1\}^n,\mu : \{(a_1,1),\dots,(a_n,n)\}\rightarrow
\Sigma,(T_d)_{d\geq 1})\] is a well-defined higher dimensional
transition system. \ep

Note that for $n = 0$, $C_0[]$, also denoted by $C_0$, is nothing else
but the higher dimensional transition system $(\{()\},\mu:\varnothing
\rightarrow \Sigma,\varnothing)$.

\bpf There is nothing to prove for $n=0,1$. So one can suppose that
$n\geq 2$. We use the characterization of Proposition~\ref{eq_hdts}.
CSA1 and the Multiset axiom are obviously satisfied. Let
\[((\epsilon_1,\dots,\epsilon_n), (a_{i_1},i_1),\dots,(a_{i_m},i_m),
(\epsilon'_1,\dots,\epsilon'_n))\] be a transition of
$C_n[a_1,\dots,a_n]$. By construction of $C_n[a_1,\dots,a_n]$, the
unique state \[(\epsilon''_1,\dots,\epsilon''_n) \in [n]\] such that
the $(p+2)$-tuple \[((\epsilon_1,\dots,\epsilon_n),
(a_{i_1},i_1),\dots,(a_{i_p},i_p),
(\epsilon''_1,\dots,\epsilon''_n))\] and the
$(m-p+2)$-tuple \[((\epsilon''_1,\dots,\epsilon''_n),
(a_{i_{p+1}},i_{p+1}),\dots,(a_{i_m},i_m),
(\epsilon'_1,\dots,\epsilon'_n))\] are transitions of
$C_n[a_1,\dots,a_n]$ is the one satisfying $\epsilon_i\leq
\epsilon''_i \leq \epsilon'_i$ for all $i\in \{1,\dots,n\}$ and
$\epsilon_i \neq \epsilon''_i$ if and only if $i\in
\{i_1,\dots,i_p\}$. So the Unique intermediate state axiom is
satisfied.  The Coherence axiom can be checked in a similar way.  \epf

Note that for every permutation $\sigma$ of $\{1,\dots,n\}$, one has
the isomorphism of weak higher dimensional transition systems
$C_n[a_1,\dots,a_n] \iso C_n[a_{\sigma(1)},\dots,a_{\sigma(n)}]$. We
must introduce $n$ distinct actions $(a_1,1),\dots,(a_n,n)$ as in
\cite{MR1461821} otherwise an object like $C_2[a,a]$ would not satisfy
the Unique intermediate state axiom.

\begin{nota} For $n\geq 1$, let $0_n = (0,\dots,0)$ ($n$-times) and
  $1_n = (1,\dots,1)$ ($n$-times). By convention, let 
  $0_0=1_0=()$. \end{nota}

\begin{nota} For $n\geq 0$, let $C_n[a_1,\dots,a_n]^{ext}$ be the weak
  higher dimensional transition system with set of states
  $\{0_n,1_n\}$, with set of actions $\{(a_1,1), \dots, (a_n,n)\}$ and
  with transitions the $(n+2)$-tuples $(0_n,(a_{\sigma(1)},\sigma(1)),
  \dots, (a_{\sigma(n)},\sigma(n)),1_n)$ for $\sigma$ running over the
  set of permutations of the set $\{1,\dots ,n\}$.
\end{nota}

\bp \label{construction_map} Let $n\geq 0$ and $a_1,\dots,a_n\in
\Sigma$. Let $X = (S,\mu:L \rightarrow \Sigma, T = \bigcup_{n\geq 1}
T_n)$ be a weak higher dimensional transition system. Let $f_0:
\{0,1\}^n \rightarrow S$ and $\widetilde{f} :
\{(a_1,1),\dots,(a_n,n)\} \linebreak[4]\rightarrow L$ be two set
maps. Then the following conditions are equivalent:
\begin{enumerate}
\item The pair $(f_0,\widetilde{f})$ induces a map of weak higher
  dimensional transition systems from $C_n[a_1,\dots,a_n]$ to $X$.
\item For every transition $((\epsilon_1, \dots, \epsilon_n),
  (a_{i_1},i_1), \dots, (a_{i_r},i_r), (\epsilon'_1, \dots,
  \epsilon'_n))$ of $C_n[a_1,\dots,a_n]$ with $(\epsilon_1, \dots,
  \epsilon_n) = 0_n$ or $(\epsilon'_1, \dots, \epsilon'_n) = 1_n$, the
  tuple \[(f_0(\epsilon_1, \dots, \epsilon_n),
  \widetilde{f}(a_{i_1},i_1), \dots, \widetilde{f}(a_{i_r},i_r),
  f_0(\epsilon'_1, \dots, \epsilon'_n))\] is a transition of $X$.
\end{enumerate} \ep

Note that the Coherence axiom plays a crucial role in the proof.

\bpf The implication $(1) \Rightarrow (2)$ is obvious. Suppose that
$(2)$ holds. Let \[((\epsilon_1, \dots, \epsilon_n),
(a_{i_{r+1}},i_{r+1}), \dots, (a_{i_{r+s}},i_{r+s}), (\epsilon'_1,
\dots, \epsilon'_n))\] be a transition of $C_n[a_1,\dots,a_n]$ with
$(\epsilon_1, \dots, \epsilon_n)\in [n]\backslash\{0_n\}$ and
$(\epsilon'_1, \dots, \epsilon'_n)\in [n]\backslash\{1_n\}$. There
exists a transition $(0_n,(a_{i_{1}},i_{1}), \dots, (a_{i_{r}},i_{r}),
(\epsilon_1, \dots, \epsilon_n))$ in $C_n[a_1,\dots,a_n]$ from $0_n$
to $(\epsilon_1, \dots, \epsilon_n)$. And there exists a transition
$((\epsilon'_1, \dots, \epsilon'_n), (a_{i_{r+s+1}},i_{r+s+1}), \dots,
(a_{i_{n}},i_{n}),1_n)$ from $(\epsilon'_1, \dots, \epsilon'_n)$ to
$1_n$ in $C_n[a_1,\dots,a_n]$. Then by construction of
$C_n[a_1,\dots,a_n]$, the two tuples $(0_n,(a_{i_{1}},i_{1}), \dots,
(a_{i_{r+s}},i_{r+s}), (\epsilon'_1, \dots, \epsilon'_n))$ and
$((\epsilon_1, \dots, \epsilon_n), (a_{i_{r+1}},i_{r+1}), \dots,
(a_{i_{n}},i_{n}),1_n)$ are two transitions of $C_n[a_1,\dots,a_n]$ as
well. Thus, the transition \[((\epsilon_1, \dots, \epsilon_n),
(a_{i_{r+1}},i_{r+1}), \dots, (a_{i_{r+s}},i_{r+s}), (\epsilon'_1,
\dots, \epsilon'_n))\] is in the closure in $\bigcup_{d\geq
  1}\{0,1\}^n\p \{(a_1,1),\dots,(a_n,n)\}^d\p \{0,1\}^n$ under the
Coherence axiom of the subset of transitions of $C_n[a_1,\dots,a_n]$
of the form $(0_n, (a_{i_1},i_1), \dots, (a_{i_r},i_r), (\epsilon'_1,
\dots, \epsilon'_n))$ or $((\epsilon_1, \dots, \epsilon_n),
(a_{i_1},i_1), \dots, (a_{i_r},i_r), 1_n)$ with $(\epsilon_1, \dots,
\epsilon_n), (\epsilon'_1, \dots, \epsilon'_n)\in [n]$. Hence, one
obtains $(2) \Rightarrow (1)$.  \epf

\bth \label{important} A weak higher dimensional transition system
satisfies the Unique intermediate state axiom if and only if it is
orthogonal to the set of inclusions
\[\{C_n[a_1,\dots,a_n]^{ext} \subset C_n[a_1,\dots,a_n], n\geq 0
\hbox{ and }a_1,\dots,a_n\in\Sigma\}.\]
\eth

\bpf \underline{Only if part}. Let $X = (S,\mu:L \rightarrow \Sigma, T
= \bigcup_{n\geq 1} T_n)$ be a weak higher dimensional transition
system satisfying the Unique intermediate state axiom.  Let $n\geq 0$
and $a_1,\dots,a_n \in \Sigma$. We have to prove that the inclusion of
weak higher dimensional transition systems $C_n[a_1,\dots,a_n]^{ext}
\subset C_n[a_1,\dots,a_n]$ induces a bijection
\[\whdts(C_n[a_1,\dots,a_n],X) \stackrel{\iso} \longrightarrow
\whdts(C_n[a_1,\dots,a_n]^{ext},X).\] This fact is trivial for $n=0$
and $n=1$ since the inclusion $C_n[a_1,\dots,a_n]^{ext} \subset
C_n[a_1,\dots,a_n]$ is an equality. Suppose now that $n\geq 2$.  Let
$f,g \in \whdts(C_n[a_1,\dots,a_n],X)$ having the same restriction to
$C_n[a_1,\dots,a_n]^{ext}$. So there is the equality $\widetilde{f} =
\widetilde{g} : \{(a_1,1),\dots,(a_n,n)\} \rightarrow L$ as set
map. Moreover, one has $f_0(0_n) = g_0(0_n)$ and $f_0(1_n) =
g_0(1_n)$. Let $(\epsilon_1,\dots,\epsilon_n) \in [n]$ be a state of
$C_n[a_1,\dots,a_n]$ different from $0_n$ and $1_n$.  Then there exist
(at least) two transitions $(0_n,(a_{i_1},i_1),\dots,(a_{i_r},i_r),
(\epsilon_1,\dots,\epsilon_n))$ and
$((\epsilon_1,\dots,\epsilon_n),(a_{i_{r+1}},i_{r+1}),\dots,(a_{i_{r+s}},i_{r+s}),
1_n)$ of $C_n[a_1,\dots,a_n]$ with $r,s\geq 1$. So the four tuples
\[(f_0(0_n),\widetilde{f}(a_{i_1},i_1),\dots,\widetilde{f}(a_{i_r},i_r),
f_0(\epsilon_1,\dots,\epsilon_n)),\]
\[(f_0(\epsilon_1,\dots,\epsilon_n),\widetilde{f}(a_{i_{r+1}},i_{r+1}),\dots,\widetilde{f}(a_{i_{r+s}},i_{r+s}),
f_0(1_n)),\]
\[(g_0(0_n),\widetilde{g}(a_{i_1},i_1),\dots,\widetilde{g}(a_{i_r},i_r),
g_0(\epsilon_1,\dots,\epsilon_n))\] and 
\[(g_0(\epsilon_1,\dots,\epsilon_n),\widetilde{g}(a_{i_{r+1}},i_{r+1}),\dots,\widetilde{g}(a_{i_{r+s}},i_{r+s}),
g_0(1_n))\] are four transitions of $X$. Since $X$ satisfies the
Unique intermediate state axiom, one obtains
$f_0(\epsilon_1,\dots,\epsilon_n) =
g_0(\epsilon_1,\dots,\epsilon_n)$. Thus $f = g$ and the set map
\[\whdts(C_n[a_1,\dots,a_n],X)  \longrightarrow
\whdts(C_n[a_1,\dots,a_n]^{ext},X)\] is one-to-one.  Let $f :
C_n[a_1,\dots,a_n]^{ext} \rightarrow X$ be a map of weak higher
dimensional transition systems. The map $f$ induces a set map
$f_0:\{0_n,1_n\} \rightarrow S$ and a set map $\widetilde{f} :
\{(a_1,1),\dots,(a_n,n)\} \rightarrow L$.  Let
$(\epsilon_1,\dots,\epsilon_n) \in [n]$ be a state of
$C_n[a_1,\dots,a_n]$ different from $0_n$ and $1_n$.  Then there exist
(at least) two transitions \[(0_n,(a_{i_1},i_1),\dots,(a_{i_r},i_r),
(\epsilon_1,\dots,\epsilon_n))\] and
\[((\epsilon_1,\dots,\epsilon_n),(a_{i_{r+1}},i_{r+1}),\dots,(a_{i_{r+s}},i_{r+s}),
1_n)\] of $C_n[a_1,\dots,a_n]$ with $r,s\geq 1$. Let us denote by
$f_0(\epsilon_1,\dots,\epsilon_n)$ the unique state of $X$ such that
\[(f_0(0_n),\widetilde{f}(a_{i_1},i_1),\dots,\widetilde{f}(a_{i_r},i_r),
f_0(\epsilon_1,\dots,\epsilon_n))\] and
\[(f_0(\epsilon_1,\dots,\epsilon_n),\widetilde{f}(a_{i_{r+1}},i_{r+1}),\dots,\widetilde{f}(a_{i_{r+s}},i_{r+s}),
f_0(1_n))\] are two transitions of $X$. Since every transition from
$0_n$ to $(\epsilon_1,\dots,\epsilon_n)$ is of the form \[(0_n,
(a_{i_{\sigma(1)}},i_{\sigma(1)}), \dots,
(a_{i_{\sigma(r)}},i_{\sigma(r)}), (\epsilon_1,\dots,\epsilon_n))\]
where $\sigma$ is a permutation of $\{1,\dots,r\}$ and since every
transition from $(\epsilon_1,\dots,\epsilon_n)$ to $1_n$ is of the
form \[((\epsilon_1,\dots,\epsilon_n),
(a_{i_{\sigma'(r+1)}},i_{\sigma'(r+1)}), \dots,
(a_{i_{\sigma'(r+s)}},i_{\sigma'(r+s)}), 1_n)\] where $\sigma'$ is a
permutation of $\{r+1,\dots,r+s\}$, one obtains a well-defined set map
$f_0 : [n] \rightarrow S$.  The pair of set maps $(f_0,\widetilde{f})$
induces a well-defined map of weak higher dimensional transition
systems by Proposition~\ref{construction_map}. Therefore the set map
\[\whdts(C_n[a_1,\dots,a_n],X)  \longrightarrow
\whdts(C_n[a_1,\dots,a_n]^{ext},X)\] is onto.

\underline{If part}. Conversely, let $X = (S,\mu:L \rightarrow \Sigma,
T = \bigcup_{n\geq 1} T_n)$ be a weak higher dimensional transition
system orthogonal to the set of inclusions $\{C_n[a_1,\dots,a_n]^{ext}
\subset C_n[a_1,\dots,a_n], n\geq 0 \hbox{ and
}a_1,\dots,a_n\in\Sigma\}$. Let $(\alpha,u_1,\dots,u_n,\beta)$ be a
transition of $X$ with $n\geq 2$. Then there exists a (unique) map
$C_n[\mu(u_1),\dots,\mu(u_n)]^{ext} \longrightarrow X$ taking the
transition \[(0_n,(\mu(u_1),1),\dots,(\mu(u_n),n),1_n)\] to
the transition $(\alpha,u_1,\dots,u_n,\beta)$. By hypothesis, this map
factors uniquely as a composite \[C_n[\mu(u_1),\dots,\mu(u_n)]^{ext}
\subset C_n[\mu(u_1),\dots,\mu(u_n)] \stackrel{g}\longrightarrow X.\]
Let $1\leq p<n$. There exists a (unique) state $\nu$ of
$C_n[\mu(u_1),\dots,\mu(u_n)]$ such that the tuples $(0_n,
(\mu(u_1),1), \dots, (\mu(u_p),p), \nu)$ and $(\nu,
(\mu(u_{p+1}),p+1), \dots, (\mu(u_n),n), 1_n)$ are two transitions of
the higher dimensional transition system $C_n[\mu(u_1), \dots,
\mu(u_n)]$ by Proposition~\ref{cas_cube}. Hence the existence of a
state $\nu_1 = g_0(\nu)$ of $X$ such that the tuples
$(\alpha,u_1,\dots,u_p,\nu_1)$ and $(\nu_1,u_{p+1},\dots,u_n,\beta)$
are two transitions of $X$. Suppose that $\nu_2$ is another state of
$X$ such that $(\alpha,u_1,\dots,u_p,\nu_2)$ and
$(\nu_2,u_{p+1},\dots,u_n,\beta)$ are two transitions of $X$. Let
$\widetilde{h} = \widetilde{g} : \{(\mu(u_1),1), \dots, (\mu(u_n),n)\}
\longrightarrow L$ be defined by $\widetilde{h}(\mu(u_i),i) =
u_i$. Let $h_0 : [n] \rightarrow S$ be defined by $h_0(\nu') =
g_0(\nu')$ for $\nu' \neq \nu$ and $h_0(\nu) = \nu_2$ (instead of
$\nu_1$). By Proposition~\ref{construction_map}, the pair of set maps
$(h_0,\widetilde{h})$ yields a well-defined map of weak higher
dimensional transition systems $h : C_n[\mu(u_1), \dots, \mu(u_n)]
\longrightarrow X$. So by orthogonality, one obtains $h = g$, and
therefore $\nu_1 = \nu_2$. Thus, the weak higher dimensional
transition system $X$ satisfies the Unique intermediate state axiom.
\epf

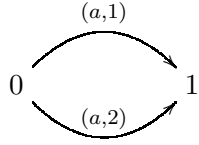
\begin{figure}
\[
\xymatrix{
0 \ar@/^20pt/[rr]^-{(a,1)} \ar@/_20pt/[rr]^-{(a,2)}&& 1 }
\]
\caption{The higher dimensional transition system $D[a]$}
\label{Da}
\end{figure}

\begin{cor} \label{ortho_hdts_class} The full subcategory $\hdts$ of
  higher dimensional transition systems is a small-orthogonality class
  of the category $\whdts$ of weak higher dimensional transition
  systems. More precisely, it is the full subcategory of objects
  orthogonal to the (unique) morphisms $D[a] \rightarrow C_1[a]$ for
  $a\in \Sigma$ and to the inclusions $C_n[a_1,\dots,a_n]^{ext}
  \subset C_n[a_1,\dots,a_n]$ for $n\geq 2$ and
  $a_1,\dots,a_n\in\Sigma$ where $D[a]$ is the higher dimensional
  transition system with set of states $\{0,1\}$, with set of labels
  $\{(a,1),(a,2)\}$, with labelling maps $\mu(a,i)=a$, and containing
  the two $1$-transitions $(0,(a,1),1)$ and $(0,(a,2),1)$ (see
  Figure~\ref{Da}).
\end{cor}

\bpf This is a consequence of Theorem~\ref{important} and
Proposition~\ref{eq_hdts}. \epf

\begin{cor} \label{stc} The full subcategory of higher dimensional
  transition systems is a full reflective locally finitely presentable
  subcategory of the category of weak higher dimensional transition
  systems. In particular, the inclusion functor $\hdts \subset \whdts$
  is limit-preserving and accessible. \end{cor}

\bpf That $\hdts$ is a full reflective locally presentable category of
$\whdts$ is a consequence of
\cite[Theorem~1.39]{MR95j:18001}. Unfortunately,
\cite[Theorem~1.39]{MR95j:18001} may be false for
$\lambda=\aleph_0$. It only enables us to conclude that the category
$\hdts$ is locally $\aleph_1$-presentable. To prove that $\hdts$ is
locally finitely presentable, we observe, thanks to
Proposition~\ref{eq_hdts}, that the notion of higher dimensional
transition system is axiomatized by the axioms of weak higher
dimensional transition system and by the two additional families of
axioms: $(\forall \alpha,u,\beta), T_{x}(\alpha,u,\beta) \Rightarrow
(\exists!  u') T_{x}(\alpha,u',\beta)$ for $x\in \Sigma$ and
\begin{multline*}(\forall \alpha,u_1,\dots,u_n,\beta),
  T_{x_1,\dots,x_n}(\alpha,u_1,\dots,u_n,\beta) \Rightarrow (\exists!
  \nu) (T_{x_1,\dots,x_p}(\alpha,u_{1},\dots,u_{p},\nu) \\\wedge
  T_{x_{p+1},\dots,x_n}(\nu,u_{p+1},\dots,u_{n},\beta))\end{multline*}
for $n\geq 2$, $1\leq p<n$ and $x_1,\dots,x_n\in\Sigma$.  So the
notion of higher dimensional transition system is axiomatized by a
limit theory, i.e. by axioms of the form $(\forall {x}), \phi({x})
\Rightarrow (\exists!  {y}\ \psi({x},{y}))$ where $\phi$ and $\psi$
are conjunctions of atomic formulas. Moreover, each symbol contains a
finite number of arguments. Hence the result as in
Theorem~\ref{loctop}.  \epf

In fact, one can easily prove that the inclusion functor $\hdts
\subset \whdts$ is finitely accessible. Let $X:I\rightarrow \hdts$ be
a directed diagram of higher dimensional transition systems.  Let $X_i
= (S_i,\mu_i:L_i \rightarrow \Sigma, T^i = \bigcup_{n\geq 1} T^i_n)$
and $X = (S,\mu:L \rightarrow \Sigma, T = \bigcup_{n\geq 1} T_n)$. The
weak higher dimensional transition system $X$ remains orthogonal to
the maps $D[a] \rightarrow C_1[a]$ for every $a\in \Sigma$ since this
property is axiomatized by the sentences $(\forall \alpha,u,\beta),
T_{x}(\alpha,u,\beta) \Rightarrow (\exists!  u')
T_{x}(\alpha,u',\beta)$ for $x\in \Sigma$. Since $\whdts$ is
topological over $\set^{\{s\}\cup \Sigma}$ by Theorem~\ref{loctop},
the colimit $\liminj X$ in $\whdts$ is the weak higher dimensional
transition system having as set of states the colimit $S = \liminj
S_i$, as set of actions the colimit $L = \liminj L_i$, as labelling
map the colimit $\mu = \liminj \mu_i$ and equipped with the final
structure of weak higher dimensional transition system. The final
structure is the set of transitions obtained by taking the closure
under the Coherence axiom of the union $\bigcup_{i} T^i$ of the image
of the $T^i$ in $\bigcup_{n\geq 1} (S \p L^n\p S)$.  Let
$(\alpha,u_1,\dots,u_n,\beta)$ be a transition of $\bigcup_{i} T^i$
with $n\geq 2$.  Let $1\leq p<n$. There exists $i\in I$ such that the
map $X_i \rightarrow \liminj X$ takes
$(\alpha^i,u^i_1,\dots,u^i_n,\beta^i)$ to
$(\alpha,u_1,\dots,u_n,\beta)$. By hypothesis, there exists a state
$\nu^i$ of $X_i$ such that $(\alpha^i,u^i_1,\dots,u^i_p,\nu^i)$ and
$(\nu^i,u^i_{p+1},\dots,u^i_n,\beta^i)$ are transitions of $X_i$. So
the map $X_i \rightarrow \liminj X$ takes $\nu^i$ to a state $\nu$ of
$\liminj X$ such that $(\alpha,u_1,\dots,u_p,\nu)$ and
$(\nu,u_{p+1},\dots,u_n,\beta)$ are transitions of $\bigcup_{i}
T^i$. Let $\nu_1$ and $\nu_2$ be two states of $\liminj X$ such that
$(\alpha,u_1,\dots,u_p,\nu_1)$, $(\nu_1,u_{p+1},\dots,u_n,\beta)$,
$(\alpha,u_1,\dots,u_p,\nu_2)$ and $(\nu_2,u_{p+1},\dots,u_n,\beta)$
are transitions of $\bigcup_{i} T^i$. Since the diagram $X$ is
directed, these four transitions come from four transitions of some
$X_j$. So $\nu_1 = \nu_2$ since $X_j$ satisfies the Unique
intermediate state axiom. Thus, the set of transitions $\bigcup_{i}
T^i$ satisfies the Unique intermediate state axiom. So by
Theorem~\ref{cube_final}, the set of transitions $\bigcup_{i} T^i$ is
the final structure and $X$ satisfies the Unique intermediate state
axiom. Therefore the inclusion functor $\hdts \subset \whdts$ is
finitely accessible.

\section{Labelled symmetric precubical sets}
\label{rappel_sym_pre_set}

The category of partially ordered sets or posets together with the
strictly increasing maps ($x<y$ implies $f(x)<f(y)$) is denoted by
$\poset$.

Let $\delta_i^\alpha : [n-1] \rightarrow [n]$ be the set map defined
for $1\leq i\leq n$ and $\alpha \in \{0,1\}$ by
$\delta_i^\alpha(\epsilon_1, \dots, \epsilon_{n-1}) = (\epsilon_1,
\dots, \epsilon_{i-1}, \alpha, \epsilon_i, \dots, \epsilon_{n-1})$.
These maps are called the \textit{face maps}.  The \textit{reduced box
  category}, denoted by $\square$, is the subcategory of $\poset$ with
the set of objects $\{[n],n\geq 0\}$ and generated by the morphisms
$\delta_i^\alpha$.  They satisfy the cocubical relations
$\delta_j^\beta \delta_i^\alpha = \delta_i^\alpha \delta_{j-1}^\beta $
for $i<j$ and for all $(\alpha,\beta)\in \{0,1\}^2$. In fact, these
algebraic relations give a presentation by generators and relations of
$\square$.

\bp \cite[Proposition~3.1]{symcub} \label{distance} Let $n\geq 1$. Let
$(\epsilon_1,\dots,\epsilon_n)$ and $(\epsilon'_1,\dots,\epsilon'_n)$
be two elements of the poset $[n]$ with $(\epsilon_1,\dots,\epsilon_n)
\leq (\epsilon'_1,\dots,\epsilon'_n)$. Then there exist
$i_1>\dots>i_{n-r}$ and $\alpha_1,\dots,\alpha_{n-r} \in \{0,1\}$ such
that $(\epsilon_1,\dots,\epsilon_n) = \delta_{i_1}^{\alpha_1} \dots
\delta_{i_{n-r}}^{\alpha_{n-r}}(0\dots 0)$ and
$(\epsilon'_1,\dots,\epsilon'_n) = \delta_{i_1}^{\alpha_1} \dots
\delta_{i_{n-r}}^{\alpha_{n-r}}(1\dots 1)$ where $r\geq 0$ is the
number of $0$ and $1$ in the arguments $0\dots 0$ and $1\dots 1$.  In
other terms, $(\epsilon_1,\dots,\epsilon_n)$ is the bottom element and
$(\epsilon'_1,\dots,\epsilon'_n)$ the top element of a $r$-dimensional
subcube of $[n]$.  \ep

\bd Let $n\geq 1$. Let $(\epsilon_1,\dots,\epsilon_n)$ and
$(\epsilon'_1,\dots,\epsilon'_n)$ be two elements of the poset
$[n]$. The integer $r$ of Proposition~\ref{distance} is called the
{\rm distance} between $(\epsilon_1,\dots,\epsilon_n)$ and
$(\epsilon'_1,\dots,\epsilon'_n)$.  Let us denote this situation by $r
= d((\epsilon_1,\dots,\epsilon_n),(\epsilon'_1,\dots,\epsilon'_n))$.
By definition, one has \[r = \sum_{i=1}^{i=n}
|\epsilon_i-\epsilon'_i|.\] \ed

\bd A set map $f:[m] \rightarrow [n]$ is {\rm adjacency-preserving} if
it is strictly increasing and if
$d((\epsilon_1,\dots,\epsilon_m),(\epsilon'_1,\dots,\epsilon'_m)) = 1$
implies
$d(f(\epsilon_1,\dots,\epsilon_m),f(\epsilon'_1,\dots,\epsilon'_m)) =
1$.  The subcategory of $\poset$ with set of objects $\{[n],n\geq 0\}$
generated by the adjacency-preserving maps is denoted by
$\widehat{\square}$. \ed

Let $\sigma_i:[n] \rightarrow [n]$ be the set map defined for $1\leq
i\leq n-1$ and $n\geq 2$ by $\sigma_i(\epsilon_1, \dots, \epsilon_{n})
= (\epsilon_1, \dots, \epsilon_{i-1},\epsilon_{i+1},\epsilon_{i},
\epsilon_{i+2},\dots,\epsilon_{n})$. These maps are called the
\textit{symmetry maps}. The face maps and the symmetry maps are
examples of adjacency-preserving maps.

\bp \cite[Proposition~A.3]{symcub} \label{carsym} Let $f:[m]
\rightarrow [n]$ be an adjacency-preserving map. The following
conditions are equivalent:
\begin{enumerate}
\item The map $f$ is a composite of face maps and symmetry maps.
\item The map $f$ is one-to-one. 
\end{enumerate}
\ep

\begin{nota} The subcategory of $\widehat{\square}$ generated by the
  one-to-one adjacency-preserving maps is denoted by $\square_S$.  In
  particular, one has the inclusions of categories 
\[\square \subset
  \square_S \subset \widehat{\square}.\] \end{nota}

By \cite[Theorem~8.1]{MR1988396}, the category $\square_S$ is the
quotient of the free category generated by the face maps
$\delta_i^\alpha$ and symmetry maps $\sigma_i$, by the following
algebraic relations:
\begin{itemize}
\item the cocubical relations $\delta_j^\beta \delta_i^\alpha =
  \delta_i^\alpha \delta_{j-1}^\beta $ for $i<j$ and for all
  $(\alpha,\beta)\in \{0,1\}^2$
\item the Moore relations for symmetry operators
  $\sigma_i\sigma_i=\id$, $\sigma_i\sigma_j\sigma_i = \sigma_j\sigma_i
  \sigma_j$ for $i=j-1$ and $\sigma_i\sigma_j=\sigma_j\sigma_i$ for
  $i<j-1$
\item the relations
  $\sigma_i\delta_j^\alpha=\delta_j^\alpha\sigma_{i-1}$ for $j<i$,
  $\sigma_i\delta_j^\alpha=\delta_{i+1}^\alpha$ for $j=i$,
  $\sigma_i\delta_j^\alpha=\delta_i^\alpha$ for $j=i+1$ and
  $\sigma_i\delta_j^\alpha=\delta_j^\alpha\sigma_{i}$ for $j>i+1$.
\end{itemize}

\bd A {\rm symmetric precubical set} is a presheaf over
$\square_S$. The corresponding category is denoted by
$\square_S^{op}\set$. If $K$ is a symmetric precubical set, then let
$K_n := K([n])$ and for every set map $f:[m] \rightarrow [n]$ of
$\square_S$, denote by $f^* : K_n \rightarrow K_m$ the corresponding
set map. \ed

Let $\square_S[n]:=\square_S(-,[n])$. It is called the
\textit{$n$-dimensional (symmetric) cube}. By the Yoneda lemma, one
has the natural bijection of sets
$\square_S^{op}\set(\square_S[n],K)\iso K_n$ for every precubical set
$K$. The \textit{boundary} of $\square_S[n]$ is the symmetric
precubical set denoted by $\de \square_S[n]$ defined by removing the
interior of $\square_S[n]$: $(\de \square_S[n])_k := (\square_S[n])_k$
for $k<n$ and $(\de \square_S[n])_k = \varnothing$ for $k\geq n$.  In
particular, one has $\de \square_S[0] = \varnothing$. An
\textit{$n$-dimensional} symmetric precubical set $K$ is a symmetric
precubical set such that $K_p = \varnothing$ for $p > n$ and $K_n \neq
\varnothing$. The labelled at most $n$-dimensional symmetric
precubical set $K_{\leq n}$ denotes the labelled symmetric precubical
set defined by $(K_{\leq n})_p = K_p$ for $p \leq n$ and $(K_{\leq
  n})_p = \varnothing$ for $p >n$.

\begin{nota} Let $f:K \rightarrow L$ be a morphism of symmetric
  precubical sets.  Let $n\geq 0$.  The set map from $K_n$ to $L_n$
  induced by $f$ will be sometimes denoted by $f_n$. \end{nota}

\begin{nota} Let $\de_i^\alpha = (\delta_i^\alpha)^*$. And let $s_i =
  (\sigma_i)^*$. \end{nota}

\bp (\cite[Proposition~A.4]{symcub}) \label{LABEL} The following data
define a symmetric precubical set denoted by $!^S\Sigma$:
\begin{itemize}
\item $(!^S\Sigma)_0=\{()\}$ (the empty word)
\item for $n\geq 1$, $(!^S\Sigma)_n=\Sigma^n$
\item $\de_i^0(a_1,\dots,a_n) = \de_i^1(a_1,\dots,a_n) =
  (a_1,\dots,\widehat{a_i},\dots,a_n)$ where the notation
  $\widehat{a_i}$ means that $a_i$ is removed.
\item $s_i(a_1,\dots,a_n) =
  (a_1,\dots,a_{i-1},a_{i+1},a_i,a_{i+2},\dots,a_n)$ for $1\leq i\leq n$.
\end{itemize}
Moreover, the symmetric precubical set $!^S\Sigma$ is orthogonal to
the set of morphisms \[\{\square_S[n] \sqcup_{\de\square_S[n]}
\square_S[n] \longrightarrow \square_S[n], n\geq 2\}.\]  \ep

\bd\label{def_lsps} A {\rm labelled symmetric precubical set (over
  $\Sigma$)} is an object of the comma category $\square_S^{op}\set
\ddownarrow !^S\Sigma$. That is, an object is a map of symmetric
precubical sets $\ell:K \rightarrow !^S\Sigma$ and a morphism is a
commutative diagram \[ \xymatrix{ K \ar@{->}[rr]\ar@{->}[rd]&& L
  \ar@{->}[ld]\\ & !^S\Sigma.&}
\]
The map $\ell$ is called the {\rm labelling map}.  The symmetric
precubical set $K$ is sometimes called the {\rm underlying symmetric
  precubical set} of the labelled symmetric precubical set. A labelled
symmetric precubical set $K \rightarrow !^S\Sigma$ is sometimes
denoted by $K$ without explicitly mentioning the labelling map. \ed

\begin{nota} Let $n\geq 1$. Let $a_1,\dots,a_n$ be labels of
  $\Sigma$. Let us denote by $\square_S[a_1,\dots,a_n] : \square_S[n]
  \rightarrow !^S\Sigma$ the labelled symmetric precubical set
  defined by \[\square_S[a_1,\dots,a_n](f) = f^*(a_1,\dots,a_n).\] And
  let us denote by $\de\square_S[a_1,\dots,a_n] : \de\square_S[n]
  \rightarrow !^S\Sigma$ the labelled symmetric precubical set
  defined as the composite \[\xymatrix{\de\square_S[a_1,\dots,a_n] :
  \de\square_S[n] \subset \square_S[n]
  \fr{\square_S[a_1,\dots,a_n]}&& !^S\Sigma.}\]
\end{nota}

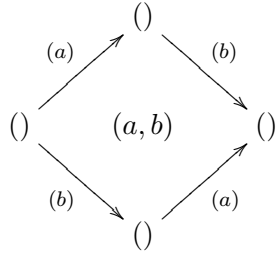
\begin{figure}
\[
\xymatrix{
& () \ar@{->}[rd]^{(b)}&\\
()\ar@{->}[ru]^{(a)}\ar@{->}[rd]_{(b)} & (a,b) & ()\\
&()\ar@{->}[ru]_{(a)}&}
\]
\caption{Concurrent execution of $a$ and $b$}
\label{concab}
\end{figure}

Figure~\ref{concab} gives the example of the labelled $2$-cube
$\square_S[a,b]$. It represents the concurrent execution of $a$ and
$b$. It is important to notice that two opposite faces of
Figure~\ref{concab} have the same label.

Since colimits are calculated objectwise for presheaves, the $n$-cubes
are finitely accessible. Since the set of cubes is a dense (and hence
strong) generator, the category of labelled symmetric precubical sets
is locally finitely presentable by \cite[Theorem~1.20 and
Proposition~1.57]{MR95j:18001}. When the set of labels $\Sigma$ is the
singleton $\{\tau\}$, the category $\square_S^{op}\set \ddownarrow
!^S\{\tau\}$ is isomorphic to the category of (unlabelled) symmetric
precubical sets since $!^S\{\tau\}$ is the terminal symmetric
precubical set.

\section{The higher dimensional automata paradigm}
\label{paradigm}

\bd A labelled symmetric precubical set $K$ satisfies {\rm the
  paradigm of higher dimensional automata (HDA paradigm)} if for every
$p\geq 2$, every commutative square of solid arrows (called a {\rm
  labelled $p$-shell} or {\rm labelled $p$-dimensional shell})
\[
\xymatrix{
  \de\square_S[p] \fr{}\fd{} && K \fd{}\\
  &&\\
  \square_S[p] \fr{} \ar@{-->}[rruu]^-{k} && !^S\Sigma}\] admits at
most one lift $k$ (i.e. a map $k$ making the two triangles
commutative).  \ed

By Definition~\ref{def_lsps}, a commutative square consisting of solid
arrows as in
\[
\xymatrix{
  \de\square_S[p] \fr{}\fd{} && K \fd{}\\
  &&\\
  \square_S[p] \fr{} \ar@{-->}[rruu]^-{k}&& !^S\Sigma}\] is equivalent
to a  diagram of labelled symmetric precubical sets consisting of the
solid arrows in 
\[
\xymatrix{\de\square_S[a_1,\dots,a_p] \fr{} \fd{} && (K\rightarrow !^S\Sigma) \\
  &&\\
  \square_S[a_1,\dots,a_p], \ar@{-->}[rruu]^-{k}&&}\] where
$(a_1,\dots,a_p)$ is the image of $\id_{[p]}$ under $\square_S[p]
\rightarrow !^S\Sigma$.  For the same reason, the existence of the
lift $k$ in the former diagram is equivalent to the existence of the
lift $k$ in the latter diagram.

\bp \label{ex_HDA} Let $n\geq 0$ and $a_1,\dots,a_n\in \Sigma$. The
labelled $n$-cube $\square_S[a_1,\dots,a_n]$ satisfies the HDA
paradigm. \ep

\bpf Consider a commutative diagram of solid arrows of the form
\[
\xymatrix{ \de\square_S[p] \fr{f}\fd{} && \square_S[n]
  \ar@{->}[dd]^-{\id_{[n] \mapsto (a_1,\dots,a_n)}}\\
  &&\\
  \square_S[p] \fr{} \ar@{-->}[rruu]^-{k} && !^S\Sigma}\] with $p\geq
2$. Then $f_0 = k_0$ as set map from $[p]$ to $[n]$. By the Yoneda
lemma, there is a bijection
$\square_S^{op}\set(\square_S[p],\square_S[n]) \iso
\square_S([p],[n])$ induced by the mapping $g\mapsto g_0$. So there
exists at most one such lift $k$. \epf

\bp \label{ortho} For a labelled symmetric precubical set $K
\rightarrow !^S\Sigma$, the following conditions are equivalent:
\begin{enumerate}
\item The labelled symmetric precubical set $K \rightarrow !^S\Sigma$
  satisfies the HDA paradigm.
\item The map $K \rightarrow !^S\Sigma$ satisfies the right lifting
  property with respect to the set of
  maps \[\left\{\square_S[p]\sqcup_{\de\square_S[p]}\square_S[p]\rightarrow
    \square_S[p], p\geq 2\right\}.\]
\item The map $K \rightarrow !^S\Sigma$ satisfies the right lifting
  property with respect to the set of
  maps \[\left\{\square_S[p]\sqcup_{\de\square_S[p]}\square_S[p] \rightarrow
    \square_S[p], p\geq 2\right\}\] and the lift is unique.
\item The labelled symmetric precubical set $K \rightarrow !^S\Sigma$
  is orthogonal to the set of maps of labelled symmetric precubical
  sets \[\left\{\square_S[a_1,\dots,a_p]
    \sqcup_{\de\square_S[a_1,\dots,a_p]}
    \square_S[a_1,\dots,a_p]\rightarrow\square_S[a_1,\dots,a_p], p\geq
    2 \hbox{ and } a_1, \dots, a_p \in \Sigma \right\}.\]
\end{enumerate}
\ep

\bpf The equivalence $(1) \iff (2)$ is due to the ``at most'' in the
definition of the HDA paradigm.  The equivalence $(3) \iff (4)$ is due
to the definition of a map of labelled symmetric precubical sets. The
implication $(3) \Longrightarrow (2)$ is obvious. The implication $(2)
\Longrightarrow (3)$ comes from the fact that for every symmetric
precubical set $K$, the set map
\[\square_S^{op}\set(\square_S[p],K) \rightarrow
\square_S^{op}\set(\square_S[p]\sqcup_{\de\square_S[p]}\square_S[p],K)\]
is one-to-one.  \epf

\begin{cor} \label{HDAreflective} The full subcategory, denoted by
  $\hda^\Sigma$, of $\square_S^{op}\set \ddownarrow !^S\Sigma$
  containing the objects satisfying the HDA paradigm is a full
  reflective locally presentable category of the category
  $\square_S^{op}\set \ddownarrow !^S\Sigma$ of labelled symmetric
  precubical sets.  In other terms, the inclusion functor $i_{\Sigma}:
  \hda^\Sigma \subset \square_S^{op}\set \ddownarrow !^S\Sigma$ has a
  left adjoint $\sh_{\Sigma}: \square_S^{op}\set \ddownarrow !^S\Sigma
  \rightarrow \hda^\Sigma$.
\end{cor}

When $\Sigma$ is the singleton $\{\tau\}$, the category $\hda^\Sigma$
will be simply denoted by $\hda$.

\bpf This is a corollary of Proposition~\ref{ortho} and
\cite[Theorem~1.39]{MR95j:18001}.  \epf

In fact the category $\hda^\Sigma$ is locally finitely presentable;
indeed, the labelled $n$-cubes for $n\geq 0$ are in $\hda^\Sigma$ by
Proposition~\ref{ex_HDA}, and one can prove that they form a dense set
of generators.

\begin{nota} When $\Sigma$ is the singleton $\{\tau\}$, let
  $i:=i_\Sigma$ and $\sh:= \sh_\Sigma$. \end{nota}

One has $i_\Sigma(K\rightarrow !^S \Sigma) \iso (i(K)\rightarrow i(!^S
\Sigma)) = (K\rightarrow !^S \Sigma)$ and $\sh_\Sigma(K\rightarrow !^S
\Sigma) \iso (\sh(K)\rightarrow \sh(!^S \Sigma) \iso !^S \Sigma)$
since the symmetric precubical set $!^S \Sigma$ already belongs to
$\hda$ by Proposition~\ref{LABEL}.

\section{Cubes as labelled symmetric precubical sets and as higher dimensional
  transition systems}
\label{iso_cube}

Let us denote by $\cube(\square^{op}_S\set\ddownarrow !^S\Sigma)$ the
full subcategory of that of labelled symmetric precubical sets
containing the labelled cubes $\square_S[a_1,\dots,a_n]$ with $n\geq
0$ and $a_1,\dots,a_n\in \Sigma$. Let us denote by $\cube(\whdts)$ the
full subcategory of that of weak higher dimensional transition systems
containing the labelled cubes $C_n[a_1,\dots,a_n]$ with $n\geq 0$ and
with $a_1,\dots,a_n\in\Sigma$. This section is devoted to proving that
these two small categories are isomorphic (cf. Theorem~\ref{ff}). Note
that $\cube(\whdts) \subset \hdts$ by Proposition~\ref{cas_cube}.

\begin{lem} \label{decomp_sym} Let $f:\square_S[m] \rightarrow
  \square_S[n]$ be a map of symmetric precubical sets. Then there
  exists a unique set map $\widehat{f} : \{1, \dots, n\} \rightarrow
  \{1, \dots, m\} \cup \{-\infty, +\infty\}$ such that $f(\epsilon_1,
  \dots, \epsilon_m) = (\epsilon_{\widehat{f}(1)}, \dots,
  \epsilon_{\widehat{f}(n)})$ for every $(\epsilon_1, \dots,
  \epsilon_m)\in [m]$ with the conventions $\epsilon_{-\infty} = 0$
  and $\epsilon_{+\infty} = 1$. Moreover, the restriction
  $\overline{f} : \widehat{f}^{-1}(\{1,\dots,m\}) \rightarrow
  \{1,\dots,m\}$ is a bijection. \end{lem}

By convention, and for the sequel, the set map $\widehat{f}$ will be
defined from $\{1, \dots, n\} \cup \{-\infty, +\infty\}$ to $\{1,
\dots, m\} \cup \{-\infty, +\infty\}$ by setting $\widehat{f}(-\infty)
= -\infty$ and $\widehat{f}(+\infty) = +\infty$.

\bpf If $\widehat{f}_1$ and $\widehat{f}_2$ are two solutions, then
one has
$(\epsilon_{\widehat{f}_1(1)},\dots,\epsilon_{\widehat{f}_1(n)}) =
(\epsilon_{\widehat{f}_2(1)},\dots,\epsilon_{\widehat{f}_2(n)})$ for
every $(\epsilon_1,\dots,\epsilon_m) \in [m]$. Let $i\in \{1, \dots,
n\}$. If $\widehat{f}_1(i) = -\infty$, then
$\epsilon_{\widehat{f}_1(i)} = 0 = \epsilon_{\widehat{f}_2(i)}$ for
every $(\epsilon_1,\dots,\epsilon_m) \in [m]$. So in this case,
$\widehat{f}_1(i) = \widehat{f}_2(i)$. For the same reason, if
$\widehat{f}_1(i) = +\infty$, then $\widehat{f}_1(i) =
\widehat{f}_2(i)$. If $\widehat{f}_1(i) \in \{1,\dots,m\}$, then
$\epsilon_{\widehat{f}_1(i)} = \epsilon_{\widehat{f}_2(i)}$ for every
$(\epsilon_1,\dots,\epsilon_m) \in [m]$. So $\widehat{f}_1(i) =
\widehat{f}_2(i)$ again. Thus, one obtains $\widehat{f}_1 =
\widehat{f}_2$. Hence there is at most one such $\widehat{f}$.
Because of the algebraic relations permuting the symmetry and face
maps recalled in Section~\ref{rappel_sym_pre_set}, the set map $f_0 :
(\square_S[m])_0 = [m] \rightarrow (\square_S[n])_0 = [n]$ factors as
a composite $[m] \rightarrow [m] \rightarrow [n]$ where the left-hand
map is a composite of symmetry maps and where the right-hand map is a
composite of face maps (see also \cite{MR1988396}). So there exists a
permutation $\sigma$ of $\{1, \dots, m\}$ such that
\[f(\epsilon_1,\dots,\epsilon_m) = \delta_{i_1}^{\alpha_1} \dots
\delta_{i_{n-m}}^{\alpha_{n-m}}(\epsilon_{\sigma(1)},\dots,\epsilon_{\sigma(m)})\]
for every $(\epsilon_1,\dots,\epsilon_m)\in [m]$.  Because of the
cocubical relations satisfied by the face maps, one can suppose that
$i_1 > i_2 > \dots > i_{n-m}$. Let $j_1 < \dots < j_m$ such that
\[\{j_1, \dots, j_m\} \cup \{i_1, \dots, i_{n-m}\} = \{1, \dots, n\}.\]
So one has $f(\epsilon_1,\dots,\epsilon_m) =
(\epsilon'_1,\dots,\epsilon'_n)$ with $\epsilon'_{i_k} = \alpha_k$ for
every $k\in \{1,\dots,n-m\}$ and $\epsilon'_{j_k} =
\epsilon_{\sigma(k)}$. Therefore, the set map $\widehat{f} :
\{1,\dots,n\} \rightarrow \{1,\dots,m\} \cup \{-\infty,+\infty\}$
defined by $\widehat{f}(i_k) = -\infty$ if $\alpha_k=0$,
$\widehat{f}(i_k) = +\infty$ if $\alpha_k=1$ and $\widehat{f}(j_k) =
\sigma(k)$ is a solution.  \epf

\begin{lem} \label{decomp_func} Let $f : \square_S[a_1,\dots,a_m]
  \rightarrow \square_S[b_1,\dots,b_n]$ and $g :
  \square_S[b_1,\dots,b_n] \rightarrow \square_S[c_1,\dots,c_p]$ be
  two maps of labelled symmetric precubical sets. Then one has
  $\widehat{g\circ f} = \widehat{f} \circ \widehat{g}$ with the
  notations of Lemma~\ref{decomp_sym}. \end{lem}

\bpf The set map $\widehat{f} : \{1,\dots ,n\} \cup
\{-\infty,+\infty\} \rightarrow \{1,\dots,m\} \cup
\{-\infty,+\infty\}$ is the unique set map such that
$f_0(\epsilon_1,\dots,\epsilon_m) =
(\epsilon_{\widehat{f}(1)},\dots,\epsilon_{\widehat{f}(n)})$ for every
$(\epsilon_1,\dots,\epsilon_m) \in [m]$ with the same notations as
above and with $\widehat{f}(-\infty) = -\infty$ and
$\widehat{f}(+\infty) = +\infty$. Therefore, one obtains the equality
$g_0(f_0(\epsilon_1,\dots,\epsilon_m)) =
g_0(\epsilon_{\widehat{f}(1)},\dots,\epsilon_{\widehat{f}(n)})$ for
every $(\epsilon_1,\dots,\epsilon_m) \in [m]$.  The set map
$\widehat{g} : \{1,\dots ,p\} \cup \{-\infty,+\infty\} \rightarrow
\{1,\dots,n\}\cup \{-\infty,+\infty\}$ is the unique set map such that
$g_0(\epsilon'_1,\dots,\epsilon'_n) =
(\epsilon'_{\widehat{g}(1)},\dots,\epsilon'_{\widehat{g}(p)})$ for
every $(\epsilon'_1,\dots,\epsilon'_n) \in [n]$ with the same
notations as above and with $\widehat{g}(-\infty) = -\infty$ and
$\widehat{g}(+\infty) = +\infty$. Let $\epsilon'_i =
\epsilon_{\widehat{f}(i)}$ for $1\leq i\leq m$. If $\widehat{g}(i) \in
\{-\infty,+\infty\}$, then $\epsilon'_{\widehat{g}(i)} =
\epsilon_{\widehat{f}(\widehat{g}(i))}$ since $\widehat{f}(-\infty) =
-\infty$ and $\widehat{f}(+\infty) = +\infty$. If $\widehat{g}(i)
\notin \{-\infty,+\infty\}$, then $\epsilon'_{\widehat{g}(i)} =
\epsilon_{\widehat{f}(\widehat{g}(i))}$ by definition of the family
$\epsilon'$. So one obtains
\[g_0(f_0(\epsilon_1,\dots,\epsilon_m)) =
g_0(\epsilon_{\widehat{f}(1)},\dots,\epsilon_{\widehat{f}(n)}) =
(\epsilon_{\widehat{f}(\widehat{g}(1))},\dots,\epsilon_{\widehat{f}(\widehat{g}(p))})\]
for every $(\epsilon_1,\dots,\epsilon_m) \in [m]$. Thus by
Lemma~\ref{decomp_sym}, one obtains $\widehat{g\circ f} = \widehat{f}
\circ \widehat{g}$. 
\epf 

Let $m,n\geq 0$ and $a_1,\dots,a_m,b_1,\dots,b_n \in \Sigma$. A map of
labelled symmetric precubical sets $f : \square_S[a_1,\dots,a_m]
\rightarrow \square_S[b_1,\dots,b_n]$ gives rise to a set
map \[\boxed{f_0 : [m] = \{0,1\}^m = \square_S[a_1,\dots,a_m]_0
  \rightarrow [n] = \{0,1\}^n = \square_S[b_1,\dots,b_n]_0}\] from the
set of states of $C_m[a_1,\dots,a_m]$ to the set of states of
$C_n[b_1,\dots,b_n]$ which belongs to $\square_S([m],[n]) =
\square_S^{op}\set(\square_S[m],\square_S[n])$. By
Lemma~\ref{decomp_sym}, there exists a unique set map $\widehat{f} :
\{1,\dots,n\} \rightarrow \{1,\dots,m\} \cup \{-\infty,+\infty\}$ such
that \[\boxed{f_0(\epsilon_1,\dots,\epsilon_m) =
  (\epsilon_{\widehat{f}(1)},\dots,\epsilon_{\widehat{f}(n)})}\] for
every $(\epsilon_1,\dots,\epsilon_m)\in [m]$ with the conventions
$\epsilon_{-\infty} = 0$ and $\epsilon_{+\infty} = 1$. Moreover, the
restriction $\overline{f} : \widehat{f}^{-1}(\{1,\dots,m\})
\rightarrow \{1,\dots,m\}$ is a bijection. Since $f :
\square_S[a_1,\dots,a_m] \rightarrow \square_S[b_1,\dots,b_n]$ is
compatible with the labelling, one necessarily has $a_i =
b_{\overline{f}^{-1}(i)}$ for every $i\in \{1,\dots,m\}$. One deduces
a set map $\widetilde{f} : \{(a_1,1),\dots,(a_m,m)\} \rightarrow
\{(b_1,1),\dots,(b_n,n)\}$ from the set of actions of
$C_m[a_1,\dots,a_m]$ to the set of actions of $C_n[b_1,\dots,b_n]$ by
setting \[\boxed{\widetilde{f}(a_i,i) =
  (b_{\overline{f}^{-1}(i)},\overline{f}^{-1}(i)) =
  (a_i,\overline{f}^{-1}(i))}.\]

\begin{lem} \label{func_map} The two set maps $f_0$ and
  $\widetilde{f}$ above defined by starting from a map of labelled
  symmetric precubical sets $f : \square_S[a_1,\dots,a_m] \rightarrow
  \square_S[b_1,\dots,b_n]$ yield a map of weak higher dimensional
  transition systems $\boxed{\T(f) : C_m[a_1,\dots,a_m] \rightarrow
  C_n[b_1,\dots,b_n]}$. \end{lem}

\bpf Let $((\epsilon_1,\dots,\epsilon_m),
(a_{i_1},i_1),\dots,(a_{i_r},i_r), (\epsilon'_1,\dots,\epsilon'_m))$
be a transition of $C_m[a_1,\dots,a_m]$. One has for every $i\in
\{1,\dots,n\}$:
\begin{itemize}
\item $\epsilon_{\widehat{f}(i)} \leq \epsilon'_{\widehat{f}(i)}$ for
  every $i \in \{1,\dots,n\}$, by definition of a transition of
  $C_m[a_1,\dots,a_m]$
\item $\epsilon_{\widehat{f}(i)} = \epsilon'_{\widehat{f}(i)}$ if $i
  \in \widehat{f}^{-1}(\{-\infty,+\infty\})$
\item $\epsilon_{\widehat{f}(i)} \neq \epsilon'_{\widehat{f}(i)}$ for
  $i \in \widehat{f}^{-1}(\{1,\dots,m\})$ if and only if
  $\widehat{f}(i) = \overline{f}(i) \in \{i_1,\dots,i_r\}$, by
  definition of a transition of $C_m[a_1,\dots,a_m]$ again.
\end{itemize}
So one has $\epsilon_{\widehat{f}(i)} \neq \epsilon'_{\widehat{f}(i)}$
if and only if $i = \overline{f}^{-1}(i_k)$ for some $k\in
\{1,\dots,r\}$. Thus, the $(d+2)$-tuple
\[((\epsilon_{\widehat{f}(1)},\dots,\epsilon_{\widehat{f}(n)}),
(a_{i_1},\overline{f}^{-1}(i_1)),\dots,(a_{i_r},\overline{f}^{-1}(i_r)),
(\epsilon'_{\widehat{f}(1)},\dots,\epsilon'_{\widehat{f}(n)}))\] is a
transition of the higher dimensional transition system
$C_n[b_1,\dots,b_n]$. \epf

\bp \label{cons_func} Let $\T(\square_S[a_1,\dots,a_n]) :=
C_n[a_1,\dots,a_n]$. Together with the mapping $f\mapsto \T(f)$
defined in Lemma~\ref{func_map}, one obtains a well-defined functor
from $\cube(\square^{op}_S\set\ddownarrow !^S\Sigma)$ to
$\cube(\whdts)$. \ep

\bpf The set map $\widehat{\id_{[m]}}$ is the inclusion $\{1,\dots
,m\} \subset \{1,\dots ,m\} \cup \{-\infty,+\infty\}$. So
\[\T(\id_{\square_S[a_1,\dots,a_n]}) = \id_{C_n[a_1,\dots,a_n]}.\]
Let $f : \square_S[a_1,\dots,a_m] \rightarrow
\square_S[b_1,\dots,b_n]$ and $g : \square_S[b_1,\dots,b_n]
\rightarrow \square_S[c_1,\dots,c_p]$ be two maps of labelled
symmetric precubical sets. The functoriality of the mapping $K \mapsto
K_{\leq 0}$ yields the equality $(g\circ f)_0 = g_0 \circ f_0$.  One
has $\widetilde{g}(\widetilde{f}(a_i,i)) =
\widetilde{g}(a_i,\overline{f}^{-1}(i)) =
(a_i,\overline{g}^{-1}(\overline{f}^{-1}(i)))$.  The integer
$N=\overline{g}^{-1}(\overline{f}^{-1}(i)) \in \{1,..,p\}$ satisfies
$i = \widehat{f}(\widehat{g}(N)) = \widehat{g\circ f}(N)$ by
Lemma~\ref{decomp_func}. So by Lemma~\ref{decomp_sym}, one has
$\overline{g\circ f}^{-1}(i) = N$. Thus, one obtains 
\[\widetilde{g}(\widetilde{f}(a_i,i)) = (a_i,\overline{g\circ
  f}^{-1}(i)) = \widetilde{g\circ f}(a_i,i).\] Hence the
functoriality. \epf

\bth \label{ff} The functor $\T : \cube(\square^{op}_S\set
\ddownarrow !^S\Sigma) \rightarrow \cube(\whdts)$ constructed in
Proposition~\ref{cons_func} is an isomorphism of categories.  \eth

\bpf Let us construct a functor $\T^{-1} : \cube(\whdts) \rightarrow
\cube(\square^{op}_S\set\ddownarrow !^S\Sigma)$ such that $\T\circ
\T^{-1} = \id_{ \cube(\whdts)}$ and $\T^{-1} \circ \T =
\id_{\cube(\square^{op}_S\set\ddownarrow !^S\Sigma)}$.

Let $\T^{-1}(C_n[a_1,\dots,a_n]) := \square_S[a_1,\dots,a_n]$ for
every $n\geq 0$ and every $a_1,\dots,a_n\in \Sigma$. Let $f:
C_m[a_1,\dots,a_m] \rightarrow C_n[b_1,\dots,b_n]$ be a map of weak
higher dimensional transition systems. By definition, it gives rise to
a set map $f_0 : [m] \rightarrow [n]$ between the set of states and to
a set map $\widetilde{f} : \{(a_1,1),\dots,(a_m,m)\} \rightarrow
\{(b_1,1),\dots,(b_n,n)\}$ between the set of actions. Since the map
$f: C_m[a_1,\dots,a_m] \rightarrow C_n[b_1,\dots,b_n]$ is compatible
with the labelling maps of the source and target higher dimensional
transition systems, one necessarily has $\widetilde{f}(a_{i},i) =
(a_i,\underline{f}(i))$ where $\underline{f} : \{1, \dots, m\}
\rightarrow \{1, \dots, n\}$ denotes a set map.  Since the
$(m+2)$-tuple \[(f_0(0,\dots,0), (a_1,\underline{f}(1)), \dots, (a_m,
\underline{f}(m)), f_0(1,\dots,1))\] is a transition of
$C_n[b_1,\dots,b_n]$, the map $\underline{f} : \{1, \dots, m\}
\rightarrow \{1, \dots, n\}$ is one-to-one. Let $\overline{f} :
\underline{f}(\{1, \dots, n\}) \rightarrow \{1, \dots, m\}$ be the
inverse map. Let $(\epsilon_1, \dots, \epsilon_m) < (\epsilon'_1,
\dots, \epsilon'_m)$ be two adjacent elements of $[m]$, more
precisely, $\epsilon_i = \epsilon'_i$ for all $i\in
\{1,\dots,m\}\backslash \{j\}$ and $0 = \epsilon_j < \epsilon'_j =
1$. Then the triple $((\epsilon_1, \dots, \epsilon_m), (a_j,j),
(\epsilon'_1, \dots, \epsilon'_m))$ is a $1$-transition of
$C_m[a_1,\dots,a_m]$. So the triple $(f_0(\epsilon_1, \dots,
\epsilon_m), (a_j,\underline{f}(j)), f_0(\epsilon'_1, \dots,
\epsilon'_m))$ is a $1$-transition of $C_n[b_1,\dots,b_n]$. Thus, the
$n$-tuples $f_0(\epsilon_1, \dots, \epsilon_m)$ and $f_0(\epsilon'_1,
\dots, \epsilon'_m)$ are adjacent in $[n]$, and the only difference is
the $\underline{f}(j)$-th coordinate. Thus, the mapping $f\mapsto f_0$
yields a set map \[(-)_0 :
\whdts(C_m[a_1,\dots,a_m],C_n[b_1,\dots,b_n]) \rightarrow
\widehat{\square}([m],[n]).\] The map $f_0:[m] \rightarrow [n]$
factors uniquely as a composite \[[m] \stackrel{\psi} \longrightarrow
[m] \stackrel{\phi} \longrightarrow [n]\] with $\phi \in \square$
since the image $f_0([m])$ is an $m$-subcube of $[n]$ (see
\cite[Proposition~3.1 and Proposition~3.11]{symcub}). Let $\phi =
\delta_{i_1}^{\alpha_1} \dots \delta_{i_{n-m}}^{\alpha_{n-m}}$ with
$i_1 > i_2 > \dots > i_{n-m}$.  Let $\widehat{f} : \{1, \dots, n\}
\cup \{-\infty, +\infty\} \rightarrow \{1, \dots, m\} \cup \{-\infty,
+\infty\}$ be the set map defined by the four mutually exclusive
cases:
\begin{itemize}
\item $\widehat{f}(-\infty) = -\infty$ and $\widehat{f}(+\infty) = +\infty$
\item $\widehat{f}(k) = \overline{f}(k)$ if $k \in \underline{f}(\{1,\dots,n\})$
\item $\widehat{f}(i_k) = -\infty$ if $\alpha_k=0$
\item $\widehat{f}(i_k) = +\infty$ if $\alpha_k=1$.
\end{itemize}
Then one has for every $m$-tuple $(\epsilon_1,\dots, \epsilon_m)$ of
$[m]$ the equality (since $\widehat{f}(\underline{f}(i)) =
i$) \[f_0(\epsilon_1,\dots, \epsilon_m) = (\epsilon_{\widehat{f}(1)},
\dots, \epsilon_{\widehat{f}(n)}).\] So $\psi\in \widehat{\square}$ is
one-to-one, and therefore equal to a composite of $\sigma_i$ maps by
Proposition~\ref{carsym}. Thus, one obtains $f_0 \in
\square_S([m],[n])$. The Yoneda bijection
\[\square_S([m],[n]) \iso \square_S^{op}\set(\square_S[m],
\square_S[n])\] takes $f_0$ to a map of symmetric precubical sets
$\T^{-1}(f) : \square_S[m] \rightarrow \square_S[n]$ preserving the
labelling. So $\T^{-1}(f)$ yields a map of labelled symmetric
precubical sets, still denoted by $\T^{-1}(f)$, from $\square_S[a_1,
\dots, a_m]$ to $\square_S[b_1, \dots, b_n]$ and it is clear that
$\T(\T^{-1}(f)) = f$ by construction of $\T$. The equality
$\T^{-1}(\T(f)) = f$ is due to the uniqueness of $\widehat{f}$ in
Lemma~\ref{decomp_sym}.  \epf

\section{Labelled symmetric precubical sets as weak
  higher dimensional transition systems}
\label{realization}

For the sequel, the category of small categories is denoted by
$\cat$. Let $H: I \longrightarrow \cat$ be a functor from a small
category $I$ to $\cat$. The \textit{Grothendieck construction} $I
\intop H$ is the category defined as follows~\cite{MR510404}: the
objects are the pairs $(i,a)$ where $i$ is an object of $I$ and $a$ is
an object of $H(i)$; a morphism $(i,a) \rightarrow (j,b)$ consists in
a map $\phi: i \rightarrow j$ and in a map $h : H(\phi)(a) \rightarrow
b$.

\begin{lem} \label{colim_id} (cf. \cite[Lemma~9.3]{symcub} and
  \cite[Lemma~A.1]{ccsprecub}) Let $I$ be a small category, and
  $i\mapsto K^i$ be a functor from $I$ to the category of labelled
  symmetric precubical sets.  Let $K = \liminj_i K^i$. Let $H: I
  \rightarrow \cat$ be the functor defined by $H(i) = {\square_S}
  \ddownarrow K^i$.  Then the functor $\iota : I\intop H \rightarrow
  {\square_S} \ddownarrow K$ defined by $\iota(i,{\square_S}[m]
  \rightarrow K^i) = ({\square_S}[m] \rightarrow K)$ is final in the
  sense of \cite{MR1712872}; that is to say the comma category $k
  \ddownarrow \iota$ is nonempty and connected for all objects $k$ of
  ${\square_S} \ddownarrow K$.  \end{lem}

\bth \label{colim1} There exists a unique colimit-preserving functor
\[\T: {\square_S}^{op}\set \ddownarrow !^S\Sigma \rightarrow \whdts\]
extending the functor $\T$ previously constructed on the full
subcategory of labelled cubes. Moreover, this functor is a left
adjoint. \eth

\bpf Let $K$ be a labelled symmetric precubical set. One necessarily
has \[\T(K) \iso \liminj_{{\square_S}[a_1,\dots,a_n] \rightarrow K}
C_n[a_1,\dots,a_n]\] hence the uniqueness. Let $K = \liminj K^i$ be a
colimit of labelled symmetric precubical sets, and denote by $I$ the
base category. By definition, one has the isomorphism
\[ \liminj \T(K^i) \iso \liminj_i
\liminj_{{\square_S}[a_1,\dots,a_n]\rightarrow K^i}
C_n[a_1,\dots,a_n].\] Consider the functor $H: I \longrightarrow \cat$
defined by $H(i) = {\square_S} \ddownarrow K^i$.  Consider the functor
$F_i: H(i) \longrightarrow \whdts$ defined by $
F_i({\square_S}[a_1,\dots,a_n] \rightarrow K^i) = C_n[a_1,\dots,a_n]$.
Consider the functor $F: I\intop H \longrightarrow \whdts$ defined by
\[F(i,{\square_S}[a_1,\dots,a_n] \rightarrow K^i) =
C_n[a_1,\dots,a_n].\] Then the composite $H(i) \subset I\intop H
\rightarrow \whdts$ is exactly $F_i$.  Therefore one has the isomorphism
\[\liminj_i \liminj_{{\square_S}[a_1,\dots,a_n]\rightarrow
  K^i} C_n[a_1,\dots,a_n] \iso
\liminj_{(i,{\square_S}[a_1,\dots,a_n]\rightarrow K^i)}
C_n[a_1,\dots,a_n]\] by \cite[Proposition
40.2]{monographie_hocolim}. The functor $\iota: I\intop H \rightarrow
{\square_S} \ddownarrow K$ defined by $\iota(i,{\square_S}[m]
\rightarrow K^i) = ({\square_S}[m] \rightarrow K)$ is final in the
sense of~\cite{MR1712872} by Lemma~\ref{colim_id}.  Therefore by
\cite[p.\ 213, Theorem 1]{MR1712872} or \cite[Theorem
14.2.5]{ref_model2}, one has the isomorphism
\[\liminj_{(i,{\square_S}[a_1,\dots,a_n]\rightarrow K^i)}
C_n[a_1,\dots,a_n] \iso \liminj_{{\square_S}[a_1,\dots,a_n]\rightarrow
  K} C_n[a_1,\dots,a_n] \iso \T(K).\] Hence the functor $\T$ is
colimit-preserving, hence the existence.  

Since the category $\square_S^{op}\set \ddownarrow !^S\Sigma$ is
locally presentable, it is co-wellpowered by
\cite[Theorem~1.58]{MR95j:18001}, and also cocomplete. The set of
labelled $n$-cubes $\{\square_S[a_1,\dots,a_n],a_1,\dots,a_n\in
\Sigma\}$ is a set of generators. So by $\hbox{SAFT}^{op}$
\cite[Corollary p126]{MR1712872}, it is a left adjoint. \epf

\bp \label{ecrasement_invisible} Let $n\geq 2$ and $a_1,\dots,a_n\in
\Sigma$. The map of labelled symmetric precubical sets
${\square_S}[a_1,\dots,a_n] \sqcup_{\de{\square_S}[a_1,\dots,a_n]}
{\square_S}[a_1,\dots,a_n] \rightarrow {\square_S}[a_1,\dots,a_n]$
induces an isomorphism of weak higher dimensional transition systems
\[\T({\square_S}[a_1,\dots,a_n] \sqcup_{\de{\square_S}[a_1,\dots,a_n]}
{\square_S}[a_1,\dots,a_n]) \iso \T({\square_S}[a_1,\dots,a_n]).\]
\ep 

\bpf Since $\T$ is colimit-preserving, one has the pushout diagram of
weak higher dimensional transition systems
\[
\xymatrix{ \T(\de{\square_S}[a_1,\dots,a_n]) \fr{}\fd{} &&
  \T({\square_S}[a_1,\dots,a_n])
  \fd{} \\
  &&\\
  \T({\square_S}[a_1,\dots,a_n]) \fr{} && \cocartesien
  \T({\square_S}[a_1,\dots,a_n] \sqcup_{\de{\square_S}[a_1,\dots,a_n]}
  {\square_S}[a_1,\dots,a_n]).}
\] 
Since $\whdts$ is topological over $\set^{\{s\}\cup \Sigma}$ by
Theorem~\ref{loctop}, the weak higher dimensional transition
system \[\T({\square_S}[a_1,\dots,a_n]
\sqcup_{\de{\square_S}[a_1,\dots,a_n]} {\square_S}[a_1,\dots,a_n])\]
is obtained by taking the colimits of the three sets of states, of the
three sets of actions and of the three labelling maps, and by endowing
the result with the final structure of weak higher dimensional
transition system. By Proposition~\ref{calcul_colim}, this final
structure is the closure under the Coherence axiom of the union of the
transitions of $\T(\de{\square_S}[a_1,\dots,a_n])$ and of the two
copies of $\T({\square_S}[a_1,\dots,a_n])$. Since the set of
transitions of $\T(\de{\square_S}[a_1,\dots,a_n])$ is included in the
set of transitions of $\T({\square_S}[a_1,\dots,a_n])$, the right-hand
vertical and bottom horizontal maps are isomorphisms.  Since the
composite
\[\T({\square_S}[a_1,\dots,a_n]) \rightarrow \T({\square_S}[a_1,\dots,a_n]
\sqcup_{\de{\square_S}[a_1,\dots,a_n]} {\square_S}[a_1,\dots,a_n])
\rightarrow \T({\square_S}[a_1,\dots,a_n])\] is the identity of
$\T({\square_S}[a_1,\dots,a_n])$, the proof is complete.  \epf

\bth \label{quotient} Let $K$ be a labelled symmetric precubical
set. The canonical map $K \rightarrow \sh_\Sigma(K)$ induces an
isomorphism of weak higher dimensional transition systems $\T(K) \iso
\T(\sh_\Sigma(K))$. \eth

\bpf By Proposition~\ref{ortho}, a labelled symmetric precubical set
$K$ belongs to $\hda^{\Sigma}$ if and only if the map $K \rightarrow
!^S\Sigma$ satisfies the right lifting property with respect to the
set of maps
\[\{{\square_S}[n] \sqcup_{\de{\square_S}[n]} {\square_S}[n]
\longrightarrow {\square_S}[n],n\geq 2\}.\] So the labelled symmetric
precubical set $\sh_\Sigma(K)$ can be obtained by a small object
argument by factoring the map $K\rightarrow !^S\Sigma$ as a composite
$K\rightarrow \sh(K) \rightarrow !^S\Sigma$ where $K \rightarrow
\sh(K)$ is a relative $\{{\square_S}[n] \sqcup_{\de{\square_S}[n]}
{\square_S}[n] \longrightarrow {\square_S}[n],n\geq 2\}$-cell complex
and where the map $\sh(K) \rightarrow !^S\Sigma$ satisfies the right
lifting property with respect to the same set of morphisms.  The small
object argument is possible by \cite[Proposition~1.3]{MR1780498} since
the category of symmetric precubical sets is locally finitely
presentable. Thanks to Proposition~\ref{ecrasement_invisible}, the
proof is complete.  \epf

\bth \label{colim2} The functor $\T:\square_S^{op}\set \ddownarrow
!^S\Sigma \rightarrow \whdts$ factors uniquely (up to isomorphism of
functors) as a composite
\[\square_S^{op}\set \ddownarrow
!^S\Sigma \stackrel{\sh_\Sigma} \longrightarrow \hda^{\Sigma}
\stackrel{\overline{\T}}\longrightarrow \whdts.\] Moreover, the
functor $\overline{\T}$ is a left adjoint. \eth

\bpf Let $\overline{\T}_1$ and $\overline{\T}_2$ be two
solutions. Then there is the isomorphism of functors
$\overline{\T}_1\circ \sh_\Sigma \iso \overline{\T}_2\circ
\sh_\Sigma$. So there are the isomorphisms of functors
$\overline{\T}_1 \iso\overline{\T}_1\circ \sh_\Sigma \circ
i_\Sigma\iso \overline{\T}_2\circ \sh_\Sigma\circ i_\Sigma\iso
\overline{\T}_2$. Let $\overline{\T} := \T \circ i_\Sigma$. Then there
is the isomorphism of functors $\overline{\T}\circ \sh_\Sigma = \T
\circ i_\Sigma \circ \sh_\Sigma \iso \T$ thanks to
Theorem~\ref{quotient}.  Hence the existence. Let $K=\liminj K_i$ be a
colimit in $\hda^\Sigma$. Then one has the sequence of natural
isomorphisms
\begin{align*}
  \overline{\T}(\liminj K_i) & \iso \overline{\T}(\liminj
  \sh_\Sigma(i_\Sigma(K_i))) & \hbox{since $K_i\iso \sh_\Sigma(i_\Sigma(K_i))$}\\
  &\iso \overline{\T}(\sh_\Sigma(\liminj i_\Sigma(K_i)))&
  \hbox{since $\sh_\Sigma$ is a left adjoint}\\
  &\iso \T(i_\Sigma(\sh_\Sigma(\liminj i_\Sigma(K_i))))& \hbox{by definition of
    $\overline{\T}$}
  \\
  &\iso \T(\liminj i_\Sigma(K_i))& \hbox{by Theorem~\ref{quotient}}\\
  &\iso \liminj \T(i_\Sigma(K_i))& \hbox{since $\T$ is colimit-preserving}\\
  &\iso \liminj \overline{\T}(K_i)&\hbox{by definition of
    $\overline{\T}$.}
\end{align*}
So the functor $\overline{\T}$ is colimit-preserving. Since the
category $\square_S^{op}\set \ddownarrow !^S\Sigma$ is locally
presentable, the functor $\overline{\T}$ is a left adjoint for the
same reason as in the proof of Theorem~\ref{colim1}.  \epf

\bd A labelled symmetric precubical set $K$ is {\rm strong} if the
weak higher dimensional transition system $\T(K)$ satisfies the Unique
intermediate state axiom. \ed

Note that a labelled symmetric precubical set $K$ is strong if and
only if $\sh_\Sigma(K)$ is strong, by Theorem~\ref{quotient}.

\bp \label{not-strong} There exists a labelled symmetric precubical
set satisfying the HDA paradigm $K$ which is not strong.  \ep

\bpf[Sketch of proof] Consider the following $1$-dimensional
(symmetric) precubical set:
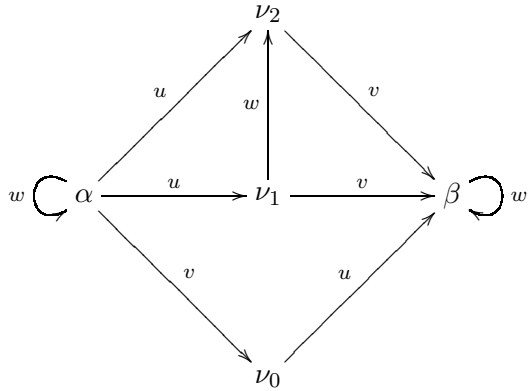
\[
\xymatrix{
  && \nu_2 \ar@{->}[rrdd]^-{v} && \\
  &&&& \\
  \ar@(ul,dl)_-{w} \alpha \ar@{->}[rruu]^-{u}
  \ar@{->}[rrdd]^-{v} \fr{u} && \fu{w}
  \nu_1 \fr{v} && \beta \ar@(ur,dr)^-{w}\\
  &&&& \\
  && \nu_0 \ar@{->}[rruu]^-{u} && }
\] 
And let us add three squares corresponding to the concurrent execution
of $u$ and $w$ (square $(\alpha,\alpha,\nu_1,\nu_2)$), of $v$ and $w$
(square $(\beta,\beta,\nu_1,\nu_2)$), and finally of $u$ and $v$
(square $(\alpha,\nu_0,\beta,\nu_1)$).  One obtains a $2$-dimensional
labelled symmetric precubical set $K$. The weak higher dimensional
transition system $\T(K)$ contains the $2$-transition
$(\alpha,u,v,\beta)$. And there exist two distinct states $\nu_1$ and
$\nu_2$ such that $(\alpha,u,\nu_1)$, $(\alpha,u,\nu_2)$,
$(\nu_1,v,\beta)$ and $(\nu_2,v,\beta)$ are $1$-transitions of the
weak higher dimensional transition system $\T(K)$.  \epf

Note that every weak higher dimensional transition system of the form
$\T(K)$ where $K$ is a labelled symmetric precubical set satisfies a
weak version of the Unique intermediate state axiom (called the
Intermediate state axiom):

\bp \label{weak-intermediate} Let $K$ be a labelled symmetric
precubical set. For every $n\geq 2$, every $p$ with $1\leq p<n$ and
every transition $(\alpha,u_1,\dots,u_n,\beta)$ of $\T(K)$, there
exists a (not necessarily unique) state $\nu$ such that both
$(\alpha,u_1,\dots,u_p,\nu)$ and $(\nu,u_{p+1},\dots,u_n,\beta)$ are
transitions. \ep

\bpf It suffices to prove that for every pushout diagram of labelled
symmetric precubical sets of the form
\[
\xymatrix{
\de\square_S[a_1,\dots,a_n] \fr{} \fd{} && K \fd{} \\
&& \\
\square_S[a_1,\dots,a_n] \fr{} && L.\cocartesien}
\]
with $n\geq 2$, if $\T(K)$ satisfies the Intermediate state axiom,
then $\T(L)$ does too. Since $\T$ is colimit-preserving by
Theorem~\ref{colim1}, one obtains the pushout diagram of weak higher
dimensional transition systems
\[
\xymatrix{
\T(\de\square_S[a_1,\dots,a_n]) \fr{f} \fd{} && \T(K) \fd{} \\
&& \\
\T(\square_S[a_1,\dots,a_n]) \fr{} && \T(L).\cocartesien}
\]
It then suffices to observe that for every $1\leq p<n$, there exists a
state $\nu_p$ of $\T(K)$ such that the tuples $(f_0(0_n),
\widetilde{f}(a_1,1), \dots, \widetilde{f}(a_p,p), \nu_p)$ and
$(\nu_p, \widetilde{f}(a_{p+1},p+1), \dots,
\widetilde{f}(a_n,n),f_0(1_n))$ are transitions of $\T(L)$: take
$\nu_p = f_0(\nu'_p)$ where $\nu'_p$ is the unique state of
$\T(\square_S[a_1,\dots,a_n])$ such that the tuples $(0_n, (a_1,1),
\dots, (a_p,p), \nu'_p)$ and $(\nu'_p, (a_{p+1},p+1), \dots,
(a_n,n),1_n)$ are transitions of $\T(\square_S[a_1,\dots,a_n])$
(cf. Proposition~\ref{cas_cube}).  \epf

One will see that for every concurrent process $P$ of every process
algebra of any synchronization algebra, the interpretation
$\square_S\llbracket P\rrbracket$ of $P$ as labelled symmetric
precubical set is always strong. In fact, it is even always a higher
dimensional transition system since CSA1 is also satisfied.

\section{Categorical property of the realization}
\label{property_func}

\bth \label{injection} Let $K$ and $L$ be two labelled symmetric
precubical sets with $L\in \hda^\Sigma$. Then the set map
\[\xymatrix{(\square_S^{op}\set\ddownarrow !^S\Sigma)(K,L) \fr{f\mapsto
  \T(f)} && \whdts(\T(K),\T(L))}\] is one-to-one.
\eth

\bpf Let $K$ and $L$ be two labelled symmetric precubical sets with
$L\in \hda^\Sigma$. Let us consider the commutative diagram of sets
\[
\xymatrix{ (\square_S^{op}\set\ddownarrow !^S\Sigma)(K,L)
  \fr{\T}\fd{(-)_{\leq 1}} &&
  \whdts(\T(K),\T(L)) \fd{}\\
  &&\\
  \set(K_{\leq 1},L_{\leq 1}) \ar@{=}[rr] && \set(\T(K)_{\leq 1},\T(L)_{\leq
    1})}
\] 
where the left-hand vertical map is the restriction to dimension $1$
and where the right-hand vertical map is the restriction of a map of
weak higher dimensional transition systems to the underlying map
between the $1$-dimensional parts, i.e. by keeping only the
$1$-dimensional transitions. The right-hand vertical map is one-to-one
by definition of a map of weak higher dimensional transition
systems. Let $f,g : K \rightrightarrows L$ be two maps of labelled
symmetric precubical sets with $f_{\leq 1} = g_{\leq 1}$. let us prove
by induction on $n\geq 1$ that $f_{\leq n} = g_{\leq n}$. The
assertion is true for $n = 1$ by hypothesis. Let us suppose that it is
true for some $n \geq 1$. Let $x:\square_S[n+1] \rightarrow K$ be a
$(n+1)$-cube of $K$. Let $\de x : \de\square_S[n+1] \subset
\square_S[n+1] \rightarrow K$. Consider the diagram of labelled
symmetric precubical sets
\[
\xymatrix{ \square_S[n+1]\sqcup_{\de\square_S[n+1]}
  \square_S[n+1]\ar@/^20pt/[rrrr]^-{f(x)\sqcup_{f(\de x)}
    g(x)}\fr{}\fd{} &&
  K_{\leq n+1} \fr{} && L \fd{}\\
  &&&&\\
  \square_S[n+1] \ar@{-->}[rrrruu]^-{k}\ar@{->}[rrrr]&&&& !^S \Sigma.}
\] 
Since $L\in \hda^\Sigma$, there exists exactly one lift $k$. Thus,
$f(x)=g(x)$ and the induction is complete.  \epf

\begin{cor} \label{injection0} The functor $\overline{\T}: \hda^\Sigma
  \rightarrow \whdts$ is faithful. \end{cor}

\bp \label{eq} The functor $\overline{\T}: \hda^\Sigma \rightarrow
\whdts$ is not full. \ep

\bpf[Sketch of proof] Let us consider the higher dimensional
transition system $C_2[u,v]$ with set of states $\{\alpha, \beta,
\nu_0, \nu_2\}$. And the inclusion of this higher dimensional
transition system to the weak higher dimensional transition system
$X=\T(K)$ given in the proof of Proposition~\ref{not-strong}:
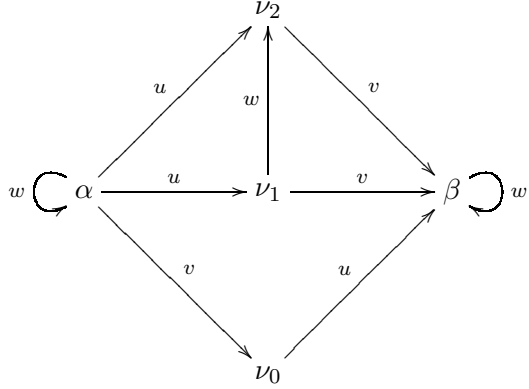
\[
\xymatrix{
  && \nu_2 \ar@{->}[rrdd]^-{v} && \\
  &&&& \\
  \ar@(ul,dl)_-{w} \alpha \ar@{->}[rruu]^-{u}
  \ar@{->}[rrdd]^-{v} \fr{u} && \fu{w}
  \nu_1 \fr{v} && \beta \ar@(ur,dr)^-{w}\\
  &&&& \\
  && \nu_0 \ar@{->}[rruu]^-{u} && }
\] 
Then this inclusion cannot come from a map of labelled symmetric
precubical sets since there are no squares in $K$ with the vertices
$\alpha,\beta,\nu_0,\nu_2$. \epf

\bth \label{bij} Let $K$ and $L$ be two labelled symmetric precubical
sets such that $L$ satisfies the HDA paradigm and such that $\T(L)$
satisfies the Unique intermediate state axiom. Then the set
map \[\xymatrix{(\square_S^{op}\set\ddownarrow !^S\Sigma)(K,L)
  \fr{f\mapsto \T(f)} && \whdts(\T(K),\T(L))}\] is bijective. \eth

\bpf First of all, let us consider the local case, i.e. when $K =
\square_S[a_1, \dots, a_m]$ is a labelled $m$-cube. It suffices to
prove that the map \[\xymatrix{(\square_S^{op}\set\ddownarrow !^S\Sigma)(K,L)
  \fr{f\mapsto \T(f)} && \whdts(\T(K),\T(L))}\] is onto since we
already know by Theorem~\ref{injection} that it is one-to-one because
$L$ satisfies the HDA paradigm. Since $\T$ is
colimit-preserving, one has the isomorphism
\[\T(L) \iso \liminj_{\square_S[b_1,\dots,b_n] \rightarrow 
  L} C_n[b_1,\dots,b_n].\] Let $f\in
\whdts(C_m[a_1,\dots,a_m],\T(L))$.
The $(m+2)$-tuple
\[(f_0(0,\dots,0), \widetilde{f}(a_1,1), \dots, \widetilde{f}(a_m,m),
f_0(1,\dots,1))\] is an $m$-transition of $\T(L)$. By
Theorem~\ref{cube_final} and since each cube $C_n[b_1,\dots,b_n]$ as
well as $\T(L)$ satisfy the Unique intermediate state axiom, there
exists a labelled cube $g:\square_S[b_1,\dots,b_n] \rightarrow L$ of
$L$ such that the $(m+2)$-tuple \[(f_0(0,\dots,0),
\widetilde{f}(a_1,1), \dots, \widetilde{f}(a_m,m), f_0(1,\dots,1))\]
comes from an $m$-transition of $C_n[b_1,\dots,b_n]$. In other terms
the composite \[C_m[a_1,\dots,a_m]^{ext} \subset C_m[a_1,\dots,a_m]
\stackrel{f} \longrightarrow \T(L)\] factors as a composite
\[C_m[a_1,\dots,a_m]^{ext} \longrightarrow C_n[b_1,\dots,b_n]
\stackrel{\T(g)}\longrightarrow \T(L).\] Since $C_n[b_1,\dots,b_n]$ is
a higher dimensional transition system by Proposition~\ref{cas_cube},
the latter map factors as a composite
\[C_m[a_1,\dots,a_m]^{ext} \subset C_m[a_1,\dots,a_m]
\stackrel{H}\longrightarrow C_n[b_1,\dots,b_n]
\stackrel{\T(g)}\longrightarrow \T(L)\] by
Theorem~\ref{important}. Since $\T(L)$ satisfies the Unique
intermediate state axiom, one obtains that the map $f\in
\whdts(C_m[a_1,\dots,a_m],\T(L))$ is equal to the composite
\[C_m[a_1,\dots,a_m] \stackrel{H}\longrightarrow C_n[b_1,\dots,b_n]
\stackrel{\T(g)} \longrightarrow \T(L)\] thanks to
Theorem~\ref{important}. By Theorem~\ref{ff}, the left-hand morphism
$H$ is of the form $\T(h)$ where $h:\square_S[a_1,\dots,a_m]
\rightarrow \square_S[b_1,\dots,b_n]$ is a map of labelled symmetric
precubical sets. Hence $f = \T(gh)$.

Let us treat now the passage from the local to the global case. Since
the functor $\T$ is colimit-preserving by Theorem~\ref{colim2}, one
has the isomorphism of weak higher dimensional transition systems
\[\T(K) \iso \liminj_{\square_S[a_1,\dots,a_m] \rightarrow 
  K} C_m[a_1,\dots,a_m].\] The set map
\begin{multline*}\limproj_{\square_S[a_1,\dots,a_m] \rightarrow 
    K} (\square_S^{op}\set\ddownarrow
  !^S\Sigma)(\square_S[a_1,\dots,a_m],L) \\\longrightarrow
  \limproj_{\square_S[a_1,\dots,a_m] \rightarrow K}
  \whdts(C_m[a_1,\dots,a_m],\T(L))\end{multline*} is bijective since
it is an inverse limit of bijections.  This completes the proof.  \epf

Note that it is also possible to prove that the set map of
Theorem~\ref{bij} is onto without using
Theorem~\ref{injection}. Indeed, the category of cubes of a labelled
symmetric precubical set is a dualizable generalized Reedy category in
the sense of \cite{gen_Reedy}. So one obtains the same result by
applying \cite[Corollary~1.7]{gen_Reedy} to the category of diagrams
from the category of cubes to the opposite $\set^{op}$ of the category
of sets and by endowing $\set$ with the unique fibrantly generated
model structure such that the fibrations are the onto maps
\cite[Theorem~4.6]{nonexistence}.

\begin{cor} Let $K$ and $L$ be two strong labelled symmetric
  precubical sets. Let us suppose that the two weak higher dimensional
  transition systems $\T(K)$ and $\T(L)$ are isomorphic. Then there is
  an isomorphism of labelled symmetric precubical sets $\sh_\Sigma(K)
  \iso \sh_\Sigma(L)$. \end{cor}

\bpf By Theorem~\ref{bij} and Theorem~\ref{quotient}, the isomorphism
$\T(\sh_\Sigma(K)) \iso \T(K) \iso \T(L) \iso \T(\sh_\Sigma(L))$ is of
the form $\T(f)$ for some map $f : \sh_\Sigma(K) \rightarrow
\sh_\Sigma(L)$ of labelled symmetric precubical sets.  And
symmetrically, there exists a map $g: \sh_\Sigma(L) \rightarrow
\sh_\Sigma(K)$ such that $\T(g)=\T(f)^{-1}$. By
Corollary~\ref{injection0}, one has $f\circ g=\id_{\sh_\Sigma(L)}$ and
$g\circ f = \id_{\sh_\Sigma(K)}$. Hence the result.  \epf

\begin{cor} \label{label_differ} Let $K$ and $L$ be two strong
  labelled symmetric precubical sets such that the weak higher
  dimensional transition systems $\T(K)$ and $\T(L)$ are
  isomorphic. Then the two weak higher dimensional transition systems
  $\T(K)$ and $\T(L)$ have the same set of actions. \end{cor}

\section{Higher dimensional transition systems are labelled symmetric
  precubical sets}
\label{coreflec}

We want to compare now the two settings of higher dimensional
transition systems and labelled symmetric precubical sets. Let us
start with some definitions and notations.
\begin{itemize}
\item $\hda^\Sigma_{hdts}$ denotes the full subcategory of
  $\hda^\Sigma$ of labelled symmetric precubical sets $K$ such that
  $\overline{\T}(K)$ is a higher dimensional transition system,
  i.e. such that the weak higher dimensional transition system
  $\overline{\T}(K)$ satisfies CSA1 and the Unique Intermediate axiom.
\item $\overline{\T}(\hda^\Sigma_{hdts})$ is the full subcategory of
  $\hdts$ of higher dimensional transition systems of the form
  $\overline{\T}(K)$ with $K\in \hda^\Sigma_{hdts}$; this subcategory
  is isomorphism-closed.
\item An action $u$ of a weak higher dimensional transition system is
  \textit{used} if there exists a transition $(\alpha,u,\beta)$.
\item The \textit{cubification} of $X\in \whdts$ is the weak higher
  dimensional transition system \[\cub(X) :=
  \liminj\limits_{C_n[a_1,\dots,a_n] \rightarrow X}
  C_n[a_1,\dots,a_n],\] the colimit being calculated in $\whdts$. Note
  that the natural map \[p_X : \cub(X) \rightarrow X\] induces a
  bijection between the set of states for any weak higher dimensional
  transition system $X$.
\end{itemize}

\bp \label{cub_prop} The cubification functor satisfies the following
properties: 
\begin{enumerate}
\item It induces a functor also denoted by $\cub$ from $\hdts$ to
  $\overline{\T}(\hda^\Sigma_{hdts})$. 
\item The natural map $p_X: \cub(X) \rightarrow X$ is an isomorphism
  for every $X \in \overline{\T}(\hda^\Sigma_{hdts})$.  
\item For every higher dimensional transition system $Y$, one has a
  natural isomorphism $\cub(\cub(Y)) \iso \cub(Y)$.
\end{enumerate}
\ep

\bpf One has, the colimit being taken in $\whdts$,
\[\cub(X) = \overline{\T}\left( \liminj\limits_{C_n[a_1,\dots,a_n]
    \rightarrow X} \square_S[a_1,\dots,a_n]\right)\] by
Theorem~\ref{ff} and Theorem~\ref{colim2}. By
Proposition~\ref{weak-intermediate}, the weak higher dimensional
transition system satisfies the Intermediate State axiom.  The
canonical map $p_X : \cub(X) \rightarrow X$ is a bijection on
states. Therefore if $X$ satisfies the Unique Intermediate State
axiom, then so does $\cub(X)$. By Theorem~\ref{cube_final}, the set of
transitions of $\cub(X)$ is the union of the transitions of the cubes
$C_n[a_1,\dots,a_n]$. So there is a bijection between the
$1$-transitions of $\cub(X)$ and the map of the form $C_1[x]
\rightarrow X$. Let $(\alpha,u,\beta)$ and $(\alpha,v,\beta)$ be two
transitions of $\cub(X)$ with $\mu(u) = \mu(v) \in \Sigma$, $\mu$
being the labelling map of $\cub(X)$. Since $X$ satisfies CSA1, one
has $\widetilde{p_X}(u) = \widetilde{p_X}(v)$. We obtain $u = v$ and
$\cub(X)$ satisfies CSA1. Hence the first assertion.

Let $K \in \hda^\Sigma_{hdts}$. Then one has 
\[\cub(\overline{\T}(K)) = \liminj\limits_{C_n[a_1,\dots,a_n] \rightarrow
  \overline{\T}(K)} C_n[a_1,\dots,a_n] \iso
\overline{\T}\left(\liminj\limits_{\square_S[a_1,\dots,a_n]
    \rightarrow K} \square_S[a_1,\dots,a_n]\right) \iso
\overline{\T}(K)\] by Theorem~\ref{ff}, Theorem~\ref{bij} and
Theorem~\ref{colim2}. Hence the second assertion.

For every higher dimensional transition system $Y$, there exists $K
\in \hda^\Sigma_{hdts}$ such that $\cub(Y) = \overline{\T}(K)$. Hence
the third assertion.  \epf

\bp \label{application0} The restriction functor $\overline{\T}:
\hda^\Sigma \rightarrow \whdts$ induces an equivalence of categories
\[\hda^\Sigma_{hdts} \simeq \overline{\T}(\hda^\Sigma_{hdts})
\simeq \hdts[\cub^{-1}]\] where $\hdts[\cub^{-1}]$ is the categorical
localization of $\hdts$ by the maps $f$ such that $\cub(f)$ is an
isomorphism. \ep

\bpf The restriction functor $\overline{\T}: \hda^\Sigma \rightarrow
\whdts$ induces an equivalence of categories $\hda^\Sigma_{hdts}
\simeq \overline{\T}(\hda^\Sigma_{hdts})$: indeed, it is faithful by
Corollary~\ref{injection0}, full by Theorem~\ref{bij} and
Proposition~\ref{eq_hdts} and essentially surjective by
construction. It remains to prove that the pair of functors \[i:
\overline{\T}(\hda^\Sigma_{hdts}) \leftrightarrows \hdts: \cub\] where
$i: \overline{\T}(\hda^\Sigma_{hdts}) \subset \hdts$ is the inclusion
functor induces an equivalence of categories
\[\overline{\T}(\hda^\Sigma_{hdts}) \simeq \hdts[\cub^{-1}].\]  That
$\cub: \hdts \rightarrow \overline{\T}(\hda^\Sigma_{hdts})$ factors
uniquely as a composite \[\hdts \rightarrow \hdts[\cub^{-1}]
\rightarrow \overline{\T}(\hda^\Sigma_{hdts})\] comes from the
universal property of the localization.  For every $X\in
\overline{\T}(\hda^\Sigma_{hdts})$, there is a natural isomorphism
$p_X : \cub(X) \iso X$ by Proposition~\ref{cub_prop}~(2). For every $Y
\in \hdts$, the map $p_Y: \cub(Y) \rightarrow Y$ is an isomorphism of
$\hdts[\cub^{-1}]$ since $\cub(p_Y)$ is an isomorphism by
Proposition~\ref{cub_prop}~(3). Hence the desired categorical
equivalence.  \epf

\bp \label{preparation0} The category
$\overline{\T}(\hda^\Sigma_{hdts})$ is a coreflective locally finitely
presentable subcategory of $\hdts$. \ep

\bpf By Proposition~\ref{cub_prop}~(2), one has the commutative
diagram of higher dimensional transition systems
\[
\xymatrix{\cub(X)\fr{\cub(f)} && \cub(Y)\fd{p_Y}\\
  &&\\
  X \fr{f} \fu{p_X^{-1}}&& Y}
\]  
for every map $f: X \rightarrow Y$ where $X$ is an object of
$\overline{\T}(\hda^\Sigma_{hdts})$ and $Y$ a higher dimensional
transition system. So the set map $h\mapsto p_Y\circ h$ from
$\hdts(X,\cub(Y))$ to $\hdts(X,Y)$ is onto. Let $f,g: X
\rightrightarrows \cub(Y)$ be two maps such that $p_Y\circ f =
p_Y\circ g$.  Since the set map $(p_Y)_0$ from the set of states of
$\cub(Y)$ to the one of $Y$ is bijective, one has $f_0 = g_0$,
i.e. $f$ and $g$ coincide on the set of states.  Let $u$ be an action
of $X$. Let $(\alpha, u ,\beta)$ be a transition of $X$: all actions
of $X$ are used since $X=\T(K)$ for some $K$. Then $(f_0(\alpha),
\widetilde{f}(u) ,f_0(\beta))$ and $(g_0(\alpha), \widetilde{g}(u)
,g_0(\beta))$ are two transitions of $\cub(Y)$. Since $f_0 = g_0$ and
since $\cub(Y)$ satisfies CSA1 by Proposition~\ref{cub_prop}~(1), one
obtains $\widetilde{f}(u) = \widetilde{g}(u)$. So $f = g$ and the map
$f\mapsto p_Y\circ f$ from $\hdts(X,\cub(Y))$ to $\hdts(X,Y)$ is
one-to-one.~\footnote{The set map $\widetilde{p_Y}$ from the set of
  actions of $\cub(Y)$ to that of $Y$ is not necessarily
  one-to-one. See Equation~(\ref{action_onto}).}

Therefore, the cubification functor $\cub: \hdts \rightarrow
\overline{\T}(\hda^\Sigma_{hdts})$ is right adjoint to the inclusion
$i: \overline{\T}(\hda^\Sigma_{hdts}) \subset \hdts$. So the category
$\overline{\T}(\hda^\Sigma_{hdts})$ is cocomplete, as a coreflective
isomorphism-closed subcategory of the cocomplete category $\hdts$.
Since $\cub(\overline{\T}(K)) = \overline{\T}(K)$, the cubes
$C_n[a_1,\dots,a_n]$ with $n\geq 0$ and $a_1,\dots, a_n \in \Sigma$
form a dense (and hence strong) generator of
$\overline{\T}(\hda^\Sigma_{hdts})$. So the category
$\overline{\T}(\hda^\Sigma_{hdts})$ is locally finitely presentable by
\cite[Theorem~1.20]{MR95j:18001}.  \epf

\bp \label{cub_hdts} A labelled symmetric precubical set $K$ belongs
to $\hda^\Sigma_{hdts}$ if and only if $K$ is orthogonal to the set of
maps \[\left\{\square_S[a_1,\dots,a_p]
  \sqcup_{\de\square_S[a_1,\dots,a_p]}
  \square_S[a_1,\dots,a_p]\rightarrow\square_S[a_1,\dots,a_p], p\geq 1
  \hbox{ and } a_1, \dots, a_p \in \Sigma \right\}\] and the weak
higher dimensional transition system $\T(K)$ satisfies the Unique
intermediate state axiom. \ep

\bpf This is a consequence of Proposition~\ref{eq_hdts}.  \epf

\bp \label{preparation} The inclusion functor $\hda^\Sigma_{hdts}
\subset \square_S^{op}\set\ddownarrow !^S\Sigma$ is limit-preserving
and finitely accessible. \ep

\bpf \underline{Limit-preserving}. Let $I$ be a small category. Let
$\underline{K}:I\rightarrow \hda^\Sigma_{hdts}$ be a diagram of
objects of $\hda^\Sigma_{hdts}$.  Then the labelled symmetric
precubical set $\limproj \underline{K}$ (limit taken in the category
of labelled symmetric precubical sets) is orthogonal to the set of
maps
\[\left\{\square_S[a_1,\dots,a_p]
  \sqcup_{\de\square_S[a_1,\dots,a_p]}
  \square_S[a_1,\dots,a_p]\rightarrow\square_S[a_1,\dots,a_p], p\geq 1
  \hbox{ and } a_1, \dots, a_p \in \Sigma \right\}\] by
\cite[Theorem~1.39]{MR95j:18001}.  It remains to prove that the weak
higher dimensional transition system $\T(\limproj \underline{K})$
satisfies the Unique intermediate state axiom by
Proposition~\ref{cub_hdts}.  Consider the canonical map of weak higher
dimensional transition systems $\T(\limproj \underline{K}) \rightarrow
\limproj (\T\circ \underline{K})$. The right-hand limit is taken in
$\hdts$ or $\whdts$ since the inclusion functor $\hdts \subset \whdts$
is a right adjoint by Corollary~\ref{stc}. Since the category $\whdts$
is topological, the set of states of $\limproj (\T\circ
\underline{K})$ is equal to the inverse limit of the sets of states of
the $\T(\underline{K}(i))$, i.e. the inverse limit of the sets of
$0$-cubes of $\underline{K}(i)$ by definition of the functor $\T$. So
the canonical map $\T(\limproj \underline{K}) \rightarrow \limproj
(\T\circ \underline{K})$ induces a bijection between the set of
states. Consequently, $\T(\limproj \underline{K})$ satisfies the
Unique Intermediate State axiom since two intermediate states for the
same transition would be mapped to the same state in $\limproj
(\T\circ \underline{K})$.  Hence, the inclusion functor
$\hda^\Sigma_{hdts} \subset \square_S^{op}\set\ddownarrow !^S\Sigma$
is limit-preserving.

\underline{Finitely accessible}. Let us now suppose that
$\underline{K}$ is directed. Then the colimit $\liminj \underline{K}$
taken in $\square_S^{op}\set\ddownarrow !^S\Sigma$ is orthogonal to
the set of maps \[\left\{\square_S[a_1,\dots,a_p]
  \sqcup_{\de\square_S[a_1,\dots,a_p]}
  \square_S[a_1,\dots,a_p]\rightarrow\square_S[a_1,\dots,a_p], p\geq 1
  \hbox{ and } a_1, \dots, a_p \in \Sigma \right\}.\] since the
inclusion functor is accessible by \cite[Theorem~1.39]{MR95j:18001}
and since every labelled cube $\square_S[a_1,\dots,a_p]$ and its
boundary $\de\square_S[a_1,\dots,a_p]$ are finitely presentable.
Moreover, one has $\T(\liminj \underline{K}) = \liminj (\T\circ
\underline{K})$ by Theorem~\ref{colim1}. So the weak higher
dimensional transition system $\T(\liminj \underline{K})$ is a higher
dimensional transition system since the inclusion functor $\hdts
\subset \whdts$ is finitely accessible as explained at the very end of
Section~\ref{small_class}. So the inclusion functor
$\hda^\Sigma_{hdts} \subset \square_S^{op}\set\ddownarrow !^S\Sigma$
is finitely accessible.  \epf

\bth \label{application1} The categorical localization
$\hdts[\cub^{-1}]$ of $\hdts$ by the maps $f$ such that $\cub(f)$ is
an isomorphism is equivalent to a full reflective locally finitely
presentable subcategory of the category of labelled symmetric
precubical sets. \eth

\bpf The theorem is a consequence of Proposition~\ref{application0},
Proposition~\ref{preparation0}, Proposition~\ref{preparation} and
\cite[Theorem~1.66]{MR95j:18001}. \epf

Let us explain what the localization $\hdts[\cub^{-1}]$ consists of.
The first effect of the cubification functor is to removed all unused
actions. Let $x \in \Sigma$. Let $\underline{x} = (\varnothing, \{x\}
\subset \Sigma, \varnothing)$ be a higher dimensional transition
system with no states and no transitions, and a unique action $x$;
then $\cub(\underline{x}) = \varnothing$. The second effect of the
cubification functor is to use different actions for two
$1$-transitions which are not related by higher dimensional cubes. For
example, one has the isomorphism
\begin{equation} \label{action_onto} C_1[x] \sqcup C_1[x] \iso \cub
  \left(\liminj\left( C_1[x] \leftarrow \underline{x} \rightarrow
      C_1[x] \right)\right).\end{equation} So, in $\hdts[\cub^{-1}]$,
two higher dimensional transition systems are isomorphic if they have
the same cubes modulo their unused actions. Given a higher dimensional
transition system $X$ all of whose actions are used, one can show that
the canonical map $\cub(X) \rightarrow X$ is bijective on states,
surjective on actions, and surjective on transitions. Using
Theorem~\ref{cube_final}, this proves that the set of transitions of a
higher dimensional transition system is always the union of the set of
transitions of its cubes.

\section{Geometric realization of a weak higher dimensional transition
system}
\label{geo}

The category $\top$ of \textit{compactly generated topological spaces}
(i.e.\ of weak Hausdorff $k$-spaces) is complete, cocomplete and
cartesian closed (more details for these kinds of topological spaces
are in \cite{MR2273730}, \cite{MR2000h:55002}, the appendix of
\cite{Ref_wH} and also in the preliminaries of \cite{model3}). For the
sequel, all topological spaces will be supposed to be compactly
generated. A \textit{compact space} is always Hausdorff.

\bd \cite{model3} A {\rm (time) flow} $X$ is a small topological
category without identity maps. The set of objects is denoted by
$X^0$.  The topological space of morphisms from $\alpha$ to $\beta$ is
denoted by $\P_{\alpha,\beta}X$. The elements of $X^0$ are also called
the {\rm states} of $X$. The elements of $\P_{\alpha,\beta}X$ are
called the {\rm (non-constant) execution paths from $\alpha$ to
  $\beta$}. A flow $X$ is {\rm loopless} if for every $\alpha\in X^0$,
the space $\P_{\alpha,\alpha}X$ is empty. \ed

\begin{nota} Let $\P X = \bigsqcup_{(\alpha,\beta)\in X^0\p X^0}
  \P_{\alpha,\beta}X$.  The topological space $\P X$ is called the
  \textit{path space} of $X$.  The source map (resp.\ the target map)
  $\P X\rightarrow X^0$ is denoted by $s$ (resp.\ $t$).
\end{nota}

\bd Let $X$ be a flow, and let $\alpha \in X^0$ be a state of $X$. The
state $\alpha$ is \textit{initial} if $\alpha\notin t(\P X)$, and the
state $\alpha$ is \textit{final} if $\alpha\notin s(\P X)$. \ed

\bd A morphism of flows $f: X \rightarrow Y$ consists of a set map
$f^0: X^0 \rightarrow Y^0$ and a continuous map $\P f: \P X
\rightarrow \P Y$ compatible with the structure. The corresponding
category is denoted by $\dtop$. \ed

The strictly associative composition law
\[
\left\{ \begin{array}{c} \P_{\alpha,\beta}X \p \P_{\beta,\gamma}X
    \longrightarrow \P_{\alpha,\gamma}X \\
    (x,y) \mapsto x*y \end{array} \right.
\]
models the composition of non-constant execution paths. The
composition law $*$ is extended in the usual way to states, that is to
constant execution paths, by $x*t(x) = x$ and $s(x)*x = x$ for every
non-constant execution path $x$.

Here are two fundamental examples of flows:
\begin{enumerate}
\item Let $S$ be a set. The flow associated with $S$, also denoted by
  $S$, has $S$ as its set of states and the empty space as its path
  space.  This construction induces a functor $\set \rightarrow \dtop$
  from the category of sets to that of flows. The flow associated with
  a set is loopless.
\item Let $(P,\leq)$ be a poset. The flow associated with $(P,\leq)$,
  also denoted by $P$, is defined as follows: the set of states of $P$
  is the underlying set of $P$; the space of morphisms from $\alpha$
  to $\beta$ is empty if $\alpha\geq \beta$ and is equal to
  $\{(\alpha,\beta)\}$ if $\alpha<\beta$, and the composition law is
  defined by $(\alpha,\beta)*(\beta,\gamma) = (\alpha,\gamma)$. This
  construction induces a functor $\poset \rightarrow \dtop$ from the
  category of posets together with the strictly increasing maps to the
  category of flows. The flow associated with a poset is loopless.
\end{enumerate}

The model structure of $\dtop$ is characterized as follows
\cite{model3}:
\begin{itemize}
\item The weak equivalences are the \textit{weak S-homotopy
    equivalences}, i.e.\ the morphisms of flows $f: X\longrightarrow
  Y$ such that $f^0: X^0\longrightarrow Y^0$ is a bijection of sets
  and such that $\P f: \P X\longrightarrow \P Y$ is a weak homotopy
  equivalence.
\item The fibrations are the morphisms of flows $f: X\longrightarrow
  Y$ such that $\P f: \P X\longrightarrow \P Y$ is a Serre
  fibration\footnote{that is, a continuous map having the RLP with
    respect to the inclusion $\mathbf{D}^n\p 0\subset \mathbf{D}^n\p
    [0,1]$ for any $n\geq 0$ where $\mathbf{D}^n$ is the
    $n$-dimensional disk.}.
\end{itemize}
This model structure is cofibrantly generated. The cofibrant
replacement functor is denoted by $(-)^{\textit{cof}}$.

A state of the flow associated with the poset
$\{\widehat{0}<\widehat{1}\}^n$ (i.e. the product of $n$ copies of
$\{\widehat{0}<\widehat{1}\}$) is denoted by an $n$-tuple of elements
of $\{\widehat{0},\widehat{1}\}$. By convention,
$\{\widehat{0}<\widehat{1}\}^0 = \{()\}$. The unique
morphism/execution path from $(x_1,\dots,x_n)$ to $(y_1,\dots,y_n)$ is
denoted by an $n$-tuple $(z_1,\dots,z_n)$ of
$\{\widehat{0},\widehat{1},*\}$ with $z_i = x_i$ if $x_i = y_i$ and
$z_i = *$ if $x_i < y_i$. For example in the flow $\{\widehat{0} <
\widehat{1}\}^2$ (cf. Figure~\ref{cube2}), one has the algebraic
relation $(*,*) = (\widehat{0},*)*(*,\widehat{1}) = (*,\widehat{0}) *
(\widehat{1},*)$.

\begin{figure}
\[
\xymatrix{ (\widehat{0},\widehat{0})
  \fr{(\widehat{0},*)}\fd{(*,\widehat{0})}
  \ar@{->}[ddrr]^-{(*,*)} && (\widehat{0},\widehat{1})\fd{(*,\widehat{1})}\\
  && \\
  (\widehat{1},\widehat{0}) \fr{(\widehat{1},*)}&&
  (\widehat{1},\widehat{1})}
\]
\caption{The flow $\{\widehat{0}<\widehat{1}\}^2$
  ($(*,*)=(\widehat{0},*)*(*,\widehat{1}) = (*,\widehat{0}) *
  (\widehat{1},*)$)}
\label{cube2}
\end{figure}

Let $\square \rightarrow \poset \subset \dtop$ be the functor defined
on objects by the mapping $[n]\mapsto \{\widehat{0}<\widehat{1}\}^n$
and on morphisms by the mapping \[\delta_i^\alpha \mapsto \lp
(\epsilon_1, \dots, \epsilon_{n-1}) \mapsto (\epsilon_1, \dots,
\epsilon_{i-1}, \alpha, \epsilon_i, \dots, \epsilon_{n-1})\rp,\] where
the $\epsilon_i$'s are elements of
$\{\widehat{0},\widehat{1},*\}$. 

Let $\square_S \rightarrow \poset \subset \dtop$ be the functor
defined on objects by the mapping $[n]\mapsto
\{\widehat{0}<\widehat{1}\}^n$ and on morphisms as follows. Let $f:[m]
\rightarrow [n]$ be a map of $\square_S$ with $m,n\geq 0$. Let
$(\epsilon_1,\dots,\epsilon_m)\in \{\widehat{0},\widehat{1},*\}^m$ be
a $r$-cube. Since $f$ is adjacency-preserving, the two elements
$f(s(\epsilon_1,\dots,\epsilon_m))$ and
$f(t(\epsilon_1,\dots,\epsilon_m))$ are respectively the initial and
final states of a unique $r$-dimensional subcube denoted by
$f(\epsilon_1,\dots,\epsilon_m)$ of $[n]$ with
$f(\epsilon_1,\dots,\epsilon_m)\in
\{\widehat{0},\widehat{1},*\}^n$. Note that the composite functor
$\square \subset \square_S \rightarrow \poset \subset \dtop$ is the
functor defined above.

\bd\cite{ccsprecub} \cite{symcub} \label{geom_rea} Let $K$ be a
labelled symmetric precubical set. The {\rm geometric realization} of
$K$ is the flow \[|K|_{flow} := \liminj_{\square_S[n] \rightarrow K}
[n]^{cof}.\] \ed

Because cubes in labelled symmetric precubical sets and in weak higher
dimensional transition systems can be identified (Theorem~\ref{ff}),
there is a well-defined functor from weak higher dimensional
transition systems to flows as follows.

\bd Let $X$ be a weak higher dimensional transition system. The {\rm
  geometric realization} of $X$ is the
flow \[|X| := \liminj_{C_n[a_1,\dots,a_n] \rightarrow X} [n]^{cof}.\] \ed

\bth \label{factor} Let $K$ be a strong labelled symmetric precubical
set satisfying the HDA paradigm, i.e. such that $\T(K)$ satisfies the
Unique intermediate state axiom. Then there is a natural isomorphism
of flows $|K|_{flow} \iso |\T(K)|$. \eth

\bpf Since $K$ is strong and since it satisfies the HDA paradigm, the
set map
\[\hda^\Sigma(\square_S[a_1,\dots,a_n],K) \rightarrow
\whdts(C_n[a_1,\dots,a_n],\T(K))\] is bijective by
Theorem~\ref{bij}. So the two colimits
\[\liminj_{\square_S[a_1,\dots,a_n]
  \rightarrow K} [n]^{cof}\] 
and 
\[\liminj_{C_n[a_1,\dots,a_n]
  \rightarrow \T(K)} [n]^{cof}\] are calculated for the same diagram of
flows.   \epf

The isomorphism $|K|_{flow} \iso |\T(K)|$ is false in
general. Consider the non-strong labelled symmetric precubical set $K$
of Proposition~\ref{not-strong}. There exists a map $C_2[u,v]
\rightarrow \T(K)$ which does not come from a square of $K$. So the
geometric realization $|\T(K)|$ contains a homotopy which is not in
$|K|_{flow}$.

\section{Process algebras and strong labelled symmetric precubical sets}
\label{proc} 

First we recall the semantics of process algebra given in
\cite{ccsprecub} and \cite{symcub}. The \textit{CCS process names} are
generated by the following syntax:
\[
P::=nil \ |\ a.P \ |\ (\nu a)P \ |\ P + P \ |\ P|| P \ |\
\rec(x)P(x)\] where $P(x)$ means a process name with one free variable
$x$. The variable $x$ must be \textit{guarded}, that is it must lie in
a prefix term $a.x$ for some $a\in\Sigma$. The set of process names is
denoted by $\proc_\Sigma$. 

The set $\Sigma\backslash\{\tau\}$, which may be empty, is supposed to
be equipped with an involution $a\mapsto \overline{a}$. In Milner's
calculus of communicating systems (CCS) \cite{0683.68008}, which is
the only case treated here, one has $a\neq \overline{a}$.  We do not
use this hypothesis. The involution on $\Sigma\backslash\{\tau\}$ is
used only in Definition~\ref{fibered} of the fibered product of two
$1$-dimensional labelled symmetric precubical sets over $\Sigma$.

\bd \label{fibered} Let $K$ and $L$ be two $1$-dimensional labelled
symmetric precubical sets.  The {\rm fibered product} of $K$ and $L$
{\rm over $\Sigma$} is the $1$-dimensional labelled symmetric
precubical set $K\p_\Sigma L$ defined as follows:
\begin{itemize}
\item $(K\p_\Sigma L)_0 = K_0\p L_0$,
\item $(K\p_\Sigma L)_1 = (K_1\p L_0) \sqcup (K_0\p L_1) \sqcup
  \{(x,y)\in K_1\p L_1, \overline{\ell(x)} = \ell(y)\}$,
\item $\de_1^\alpha(x,y) = (\de_1^\alpha(x),y)$ for every $(x,y)\in
  K_1\p L_0$,
\item $\de_1^\alpha(x,y) = (x,\de_1^\alpha(y))$ for every $(x,y)\in
  K_0\p L_1$,
\item $\de_1^\alpha(x,y) = (\de_1^\alpha(x),\de_1^\alpha(y))$ for every $(x,y)\in
  K_1\p L_1$,
\item $\ell(x,y)=\ell(x)$ for every $(x,y)\in K_1\p L_0$,
\item $\ell(x,y)=\ell(y)$ for every $(x,y)\in K_0\p L_1$,
\item $\ell(x,y)=\tau$ for every $(x,y)\in K_1\p L_1$ with
  $\overline{\ell(x)} = \ell(y)$.
\end{itemize}
The $1$-cubes $(x,y)$ of $(K\p_\Sigma L)_1\cap (K_1\p L_1)$ are called
{\rm synchronizations} of $x$ and $y$.  \ed

\begin{figure}
\[
\xymatrix{
  &\ar@{->}[rr]|{a} &&\\
  \ar@{->}[rrru]|-{\tau}\ar@{->}[ru]|-{\overline{a}}\ar@{->}[rr]&& 
\ar@{->}[ru]|-{\overline{a}}&\\
  &\ar@{->}[uu] \ar@{->}[rr]&&\ar@{->}[uu]|-{b}\\
  \ar@{->}[uu]|-{b}\ar@{->}[ru]|-{\overline{a}}\ar@{->}[rrru]|-{\tau}
  \ar@{->}[rr]|-{a}&&\ar@{->}[uu]\ar@{->}[ru]|-{\overline{a}}& }\]
\caption{Representation of $\square_S[a,b]_{\leq 1}\p_\Sigma
  \square_S[\overline{a}]$}
\label{3sync00}
\end{figure}

\bd A labelled symmetric precubical set $\ell:K\rightarrow !^S\Sigma$
{\rm decorated by process names} is a labelled precubical set together
with a set map $d:K_0 \rightarrow \proc_\Sigma$ called the {\rm
  decoration}. \ed

Let $(\square_S)_{n} \subset \square_S$ be the full subcategory of
$\square_S$ containing the $[p]$ only for $p\leq n$.  By
\cite[Proposition~5.4]{symcub}, the truncation functor
$\square_S^{op}\set \ddownarrow !^S\Sigma\rightarrow
(\square_S)_{n}^{op}\set\ddownarrow !^S\Sigma$ has a right adjoint
$\cosk_n^{\square_S,\Sigma}:(\square_S)_{n}^{op}\set\ddownarrow
!^S\Sigma \rightarrow \square_S^{op}\set\ddownarrow !^S\Sigma$.

\bd Let $K$ be a $1$-dimensional labelled symmetric precubical set
with $K_0 = [p]$ for some $p\geq 0$. The {\rm labelled symmetric
  directed coskeleton} of $K$ is the labelled symmetric precubical set
$\COSK_S^\Sigma(K)$ defined as the subobject of
$\cosk^{\square_S,\Sigma}_1(K)$ such that:
\begin{itemize}
\item $\COSK_S^\Sigma(K)_{\leq 1} = \cosk^{\square_S,\Sigma}_1(K)_{\leq 1}$
\item for every $n\geq 2$, $x\in \cosk^{\square_S,\Sigma}_1(K)_n$ is an $n$-cube
  of $\COSK_S^\Sigma(K)$ if and only if the set map
  $x_0:[n]\rightarrow [p]$ is {\rm non-twisted}, i.e. $x_0: [n]
  \rightarrow [p]$ is a composite\footnote{The factorization is
    necessarily unique.}
\[x_0: [n] \stackrel{\phi}\longrightarrow [q]
\stackrel{\psi}\longrightarrow [p],\] where $\psi$ is a morphism of
the small category $\square$ and where $\phi$ is of the form
\[(\epsilon_1,\dots,\epsilon_{n}) \mapsto
(\epsilon_{i_1},\dots,\epsilon_{i_q})\] such that
$\{1,\dots,n\}\subset \{i_1,\dots,i_q\}$.
\end{itemize} \ed

Let us recall that for every $m,n\geq 0$ and
$a_1,\dots,a_m,b_1,\dots,b_n \in \Sigma$, the labelled symmetric
precubical set $\COSK_S^{\Sigma}(\square_S[a_1,\dots,a_m]_{\leq 1}
\p_\Sigma \square_S[b_1,\dots,b_n]_{\leq 1})$ satisfies the HDA
paradigm. In particular, one has the isomorphism of labelled symmetric
precubical sets \[\boxed{\square_S[a_1,\dots,a_m] \iso
  \COSK_S^{\Sigma}(\square_S[a_1,\dots,a_m]_{\leq 1})}.\]

\bd Let $K$ and $L$ be two labelled symmetric precubical
sets. The {\rm tensor product with synchronization} (or {\rm
  synchronized tensor product}) of $K$ and $L$ is
\[K \ot_\Sigma L := \liminj_{\square_S[a_1,\dots,a_m]\rightarrow K}
\liminj_{\square_S[b_1,\dots,b_n]\rightarrow L}
\COSK_S^{\Sigma}(\square_S[a_1,\dots,a_m]_{\leq 1} \p_\Sigma
\square_S[b_1,\dots,b_n]_{\leq 1}).\] \ed

Let us define by induction on the syntax of the CCS process name $P$
the decorated labelled symmetric precubical set $\square_S\llbracket
P\rrbracket$ (see \cite{ccsprecub} for further explanations). The
labelled symmetric precubical set $\square_S\llbracket P\rrbracket$
has always a unique initial state canonically decorated by the process
name $P$ and its other states will be decorated as well in an
inductive way. Therefore for every process name $P$,
$\square_S\llbracket P\rrbracket$ is an object of the double comma
category $\{i\}\ddownarrow \square_S^{op}\set \ddownarrow !^S\Sigma$.
One has $\square_S\llbracket nil\rrbracket:=\square_S[0]$,
$\square_S\llbracket \mu.nil\rrbracket:=\mu.nil
\stackrel{(\mu)}\longrightarrow nil$, $\square_S\llbracket
P+Q\rrbracket := \square_S\llbracket P\rrbracket \oplus
\square_S\llbracket Q\rrbracket$ with the binary coproduct taken in
$\{i\}\ddownarrow \square_S^{op}\set \ddownarrow !^S\Sigma$, the
pushout diagram of symmetric precubical sets
\[\xymatrix{
    \square_S[0]=\{0\} \ar@{->}[r]^-{0\mapsto nil}
    \ar@{->}[d]^-{0\mapsto P} & 
\square_S\llbracket \mu.nil\rrbracket \ar@{->}[d] \\
    \square_S\llbracket P\rrbracket \ar@{->}[r] & \cocartesien
    {\square_S\llbracket \mu.P\rrbracket},}\]
  the pullback diagram of symmetric precubical sets
  \[\xymatrix{ \square_S\llbracket (\nu a) P\rrbracket \ar@{->}[r]
    \ar@{->}[d] \cartesien &
    \square_S\llbracket P\rrbracket \ar@{->}[d] \\
    !^S(\Sigma\backslash \{a,\overline{a}\}) \ar@{->}[r] & !^S\Sigma,}
\]
the formula giving the interpretation of the parallel composition with
synchronization \[\square_S\llbracket P||Q\rrbracket :=
\square_S\llbracket P\rrbracket \ot_\Sigma \square_S\llbracket
Q\rrbracket\] and finally $\square_S\llbracket \rec(x)P(x)\rrbracket$
defined as the least fixed point of $P(-)$.  The condition imposed on
$P(x)$ implies that for all process names $Q_1$ and $Q_2$ with
$\square_S \llbracket Q_1\rrbracket \subset \square_S \llbracket
Q_2\rrbracket$, one has $\square_S \llbracket P(Q_1)\rrbracket \subset
\square_S \llbracket P(Q_2)\rrbracket$. So by starting from the
inclusion of labelled symmetric precubical sets $\square_S\llbracket
nil\rrbracket \subset \square_S\llbracket P(nil)\rrbracket$ given by
the unique initial state of $\square_S\llbracket P(nil)\rrbracket$,
the labelled symmetric precubical set
\[\square_S\llbracket
\rec(x)P(x)\rrbracket:=\liminj\limits_n \square_S\llbracket
P^n(nil)\rrbracket \iso \bigcup_{n\geq 0} \square_S\llbracket
P^n(nil)\rrbracket\] will be equal to the least fixed point of $P(-)$.

\bp \label{c} Let $m,n\geq 0$ and $a_1,\dots,a_m,b_1,\dots,b_n\in
\Sigma$. The weak higher dimensional transition system
$\T(\COSK_S^{\Sigma}(\square_S[a_1,\dots,a_m]_{\leq 1} \p_\Sigma
\square_S[b_1,\dots,b_n]_{\leq 1}))$ is a higher dimensional
transition system. \ep

\bpf The proof is similar to the proof of Proposition~\ref{cas_cube}.
\epf

\bth \label{strong} For every CCS process name $P$, the labelled
symmetric precubical set $\square_S \llbracket P \rrbracket$ belongs
to $\hda^\Sigma$ and the weak higher dimensional transition system
$\T(\square_S \llbracket P \rrbracket)$ satisfies CSA1 and the Unique
intermediate state axiom, i.e. $\T(\square_S \llbracket P \rrbracket)
\in \hdts$. \eth

\bpf[Sketch of proof] That $\square_S \llbracket P \rrbracket$ belongs
to $\hda^\Sigma$ is proved by induction on the syntax of $P$, as in
\cite[Theorem~5.2]{ccsprecub}. If $\square_S\llbracket P\rrbracket$
and $\square_S\llbracket Q\rrbracket$ belong to $\hda^\Sigma$, then
$\square_S\llbracket P+Q\rrbracket$ belongs to $\hda^\Sigma$ since
every map $\de\square_S[a_1,\dots,a_p]\rightarrow \square_S\llbracket
P+Q\rrbracket$ with $p\geq 2$ factors as a composite
$\de\square_S[a_1,\dots,a_p]\rightarrow \square_S\llbracket
P\rrbracket \rightarrow \square_S\llbracket P+Q\rrbracket$ or as a
composite $\de\square_S[a_1,\dots,a_p]\rightarrow \square_S\llbracket
Q\rrbracket \rightarrow \square_S\llbracket P+Q\rrbracket$. If
$\square_S\llbracket P\rrbracket$ belongs to $\hda^\Sigma$, then
$\square_S\llbracket (\nu a P)\rrbracket$ belongs to $\hda^\Sigma$
since $\square_S\llbracket (\nu a) P\rrbracket \subset
\square_S\llbracket P\rrbracket$. If for every $n\geq 0$, the labelled
symmetric precubical set $\square_S\llbracket P^n(nil)\rrbracket$
belongs to $\hda^\Sigma$, then $\square_S\llbracket
\rec(x)P(x)\rrbracket$ belongs to $\hda^\Sigma$ since the inclusion
functor $\hda^\Sigma \subset \square_S^{op}\set\ddownarrow !^S\Sigma$
is accessible by Corollary~\ref{HDAreflective}. Finally, let $P$ and
$Q$ be two process names such that both $\square_S \llbracket P
\rrbracket$ and $\square_S \llbracket Q \rrbracket$ belong to
$\hda^\Sigma$. For a given map $\de\square_S[a_1,\dots,a_p]\rightarrow
\square_S \llbracket P || Q \rrbracket$ with $p\geq 2$, the category
$\de\square_S[a_1,\dots,a_p] \ddownarrow (\square_S\p \square_S)
\ddownarrow \square_S\llbracket P || Q \rrbracket$ has an initial
object~\footnote{There is an error in \cite[Theorem~5.2]{ccsprecub},
  which says that this small category is directed. It should say that
  this category always has an initial object.} otherwise $\square_S
\llbracket P \rrbracket$ or $\square_S \llbracket Q \rrbracket$ would
not satisfy the HDA paradigm. Hence the labelled symmetric precubical
set $\square_S \llbracket P || Q \rrbracket$ satisfies the HDA
paradigm too, since the labelled symmetric precubical set
$\COSK_S^{\Sigma}(\square_S[a_1,\dots,a_m]_{\leq 1} \p_\Sigma
\square_S[b_1,\dots,b_n]_{\leq 1})$ does for every $m,n\geq 0$ and for
every $a_1,\dots,a_m,b_1,\dots,b_n \in \Sigma$.
 
It is clear that CSA1 is always satisfied by $\square_S \llbracket P
\rrbracket$. That $\square_S \llbracket P \rrbracket$ is a strong
labelled symmetric precubical set, i.e. that $\T(\square_S \llbracket
P \rrbracket)$ satisfies the Unique intermediate state axiom, is
proved by induction on the syntax of $P$ as follows. It is obvious
that if $\square_S \llbracket P \rrbracket$ and $\square_S \llbracket
Q \rrbracket$ are strong, then $\square_S \llbracket P+Q \rrbracket$
is strong too.  It is also obvious that $\square_S \llbracket (\nu a)
P \rrbracket$ is strong since the weak higher dimensional transition
system $\T(\square_S \llbracket (\nu a) P \rrbracket)$ is included in
the higher dimensional transition system $\T(\square_S \llbracket P
\rrbracket)$.  If for every $n\geq 0$, the labelled symmetric
precubical set $\square_S\llbracket P^n(nil)\rrbracket$ is strong,
then $\square_S\llbracket \rec(x)P(x)\rrbracket$ is strong too by
Theorem~\ref{application1}.  It remains to prove that if the two weak
higher dimensional transition systems $\T(\square_S \llbracket P
\rrbracket)$ and $\T(\square_S \llbracket Q \rrbracket)$ satisfy the
Unique intermediate axiom, then the weak higher dimensional transition
systems $\T(\square_S \llbracket P \rrbracket \ot_\Sigma \square_S
\llbracket Q \rrbracket)$ does as well. Since $\T$ is
colimit-preserving by Theorem~\ref{colim1}, the weak higher
dimensional transition system $\T(\square_S \llbracket P \rrbracket
\ot_\Sigma \square_S \llbracket Q \rrbracket)$ is isomorphic to
\[\liminj_{\square_S[a_1,\dots,a_m]\rightarrow K}
\liminj_{\square_S[b_1,\dots,b_n]\rightarrow L}
\T(\COSK_S^{\Sigma}(\square_S[a_1,\dots,a_m]_{\leq 1} \p_\Sigma
\square_S[b_1,\dots,b_n]_{\leq 1})).\] By Theorem~\ref{cube_final} and
Proposition~\ref{c}, an $n$-transition of the higher dimensional
transition system $\T(\COSK_S^{\Sigma}(\square_S[a_1,\dots,a_m]_{\leq
  1} \p_\Sigma \square_S[b_1,\dots,b_n]_{\leq 1}))$ is of the
form \[((\alpha,\beta), (u_1,v_1), \dots, (u_n,v_n),
(\gamma,\delta))\] with three mutually exclusive cases for the
$(u_i,v_i)$: (1) Both $u_i$ and $v_i$ are actions of respectively
$\T(\square_S \llbracket P \rrbracket)$ and $\T(\square_S \llbracket Q
\rrbracket)$; in this case $u_i = \overline{v_i}$ and $\mu(u_i,v_i) =
\tau$; this case corresponds to a synchronization. (2) $u_i$ is an
action of $\T(\square_S \llbracket P \rrbracket)$ and $v_i$ is a state
of $\T(\square_S \llbracket Q \rrbracket)$. (3) $u_i$ is a state of
$\T(\square_S \llbracket P \rrbracket)$ and $v_i$ is an action of
$\T(\square_S \llbracket Q \rrbracket)$.  For such an $n$-transition,
the tuples obtained from $(\alpha, u_1, \dots, u_n, \gamma)$ and
$(\beta, v_1, \dots, v_n, \delta)$ by removing the $u_i$ and $v_i$
which are states are transitions of respectively $\T(\square_S
\llbracket P \rrbracket)$ and $\T(\square_S \llbracket Q
\rrbracket)$. So the union of the transitions of the
$\T(\COSK_S^{\Sigma}(\square_S[a_1,\dots,a_m]_{\leq 1} \p_\Sigma
\square_S[b_1,\dots,b_n]_{\leq 1}))$ satisfies the Unique intermediate
state axiom since $\T(\square_S \llbracket P \rrbracket)$ and
$\T(\square_S \llbracket Q \rrbracket)$ do. So by
Theorem~\ref{cube_final} again, this union is the final structure,
that is the colimit. Hence, the weak higher dimensional transition
system $\T(\square_S \llbracket P \rrbracket \ot_\Sigma \square_S
\llbracket Q \rrbracket)$ satisfies the Unique intermediate state
axiom.  \epf

\begin{cor} \label{final} The mapping taking each CCS process name $P$
  to the flow $|\square_S \llbracket P \rrbracket|_{flow}$ factors
  through the category of higher dimensional transition systems.
\end{cor}

\section{Concluding remarks and perspectives}

\begin{figure}
\[
\xymatrix{
\proc_\Sigma \fr{\square_S\llbracket -\rrbracket} &&
\hda^\Sigma_{hdts} \fr{\iso} \fd{\subset} &&
\hdts[\cub^{-1}]  \ar@{->}[rrdd]^-{|-|_{flow}}\fd{\subset}&& \\
&&&&&&\\
&& \square_S^{op}\set\ddownarrow !^S\Sigma \fr{\T} && \whdts \fr{|-|} && \dtop}
\] 
\caption{Recapitulation: one has the inclusions of full subcategories
  $\hdts[\opt^{-1}] \subset \hdts \subset
  \whdts$}
\label{recap}
\end{figure}

The commutative diagram of Figure~\ref{recap} summarizes the two main
results of this paper. In $\hdts[\cub^{-1}]$, two higher dimensional
transition systems are isomorphic if they have the same cubes modulo
their unused actions. This category is equivalent to a full
coreflective subcategory of the category $\hdts$ of higher dimensional
transition systems, and the latter is a reflective full subcategory of
that of weak higher dimensional transition systems $\whdts$. The
category $\hdts[\cub^{-1}]$ is also equivalent to $\hda^\Sigma_{hdts}$
which is a full reflective subcategory of that of labelled symmetric
precubical sets, and even a full reflective subcategory of those
satisfying the HDA paradigm ($\hda^\Sigma$). The inclusion
$\hda^\Sigma_{hdts} \subset \hda^\Sigma$ is strict since the
non-strong labelled symmetric precubical set $K$ used for proving
Proposition~\ref{not-strong} satisfies the HDA paradigm.

All these constructions illustrate the expressiveness of the category
of flows and of the other topological models of concurrency. Indeed,
using the geometric realization functors from $\hdts$ to $\dtop$, one
can associate a flow with any transition system with independence,
with any Petri net, with any domain of configurations of prime event
structures \cite{MR1461821}, and of course with any process algebra as
already explained in \cite{ccsprecub}.

It would be interesting to find a geometric sufficient condition for a
labelled symmetric precubical set $K$ to be strong, for example by
proving that $\hda^\Sigma_{hdts}$ is a small-orthogonality class. It
would be also interesting to find the analogue of the notion of weak
higher dimensional transition system for the labelled symmetric
transverse precubical sets (the presheaves over $\widehat{\square}$)
introduced in \cite{symcub}. Weak higher dimensional transition
systems are transition systems indexed by finite multisets of
actions. The analogous notion for labelled symmetric transverse
precubical sets should be a notion of transition system indexed by
partially ordered finite multisets of actions. By restricting to
transitions labelled by finite multisets endowed with a discrete
ordering, one should get back a weak higher dimensional transition
system. Understanding the link between labelled transverse symmetric
precubical sets and higher dimensional transition systems is necessary
since the space of morphisms of flows from $|\square_S[m]|_{flow}$ to
itself for $m\geq 0$ is homotopy equivalent to
$\widehat{\square}([m],[m])$, not to $\square_S([m],[m])$, and the
inclusion $\square_S([m],[m]) \subset \widehat{\square}([m],[m])$ is
strict for $m\geq 2$. In particular, it contains the set map
$(\epsilon_1,\dots,\epsilon_m) \mapsto
(\max(\epsilon_1,\epsilon_2),\min(\epsilon_1,\epsilon_2),\epsilon_3,
\dots,\epsilon_m)$. These questions will be hopefully the subject of
future works.


\end{document}